\documentclass[a4paper,12pt]{article}
\usepackage[
  textheight=190mm,
  left=17.5mm,
  right=17.5mm
]{geometry}

\usepackage{graphicx} 
\usepackage[english]{babel}
\usepackage[utf8]{inputenc}
\usepackage{csquotes}
\usepackage{boxedminipage}

\usepackage{mathtools}
\usepackage{amssymb}
\usepackage{amsmath}
\usepackage{amsthm}
\usepackage{bbm}
\usepackage{xurl}
\usepackage{enumerate}
\usepackage{tikz}
\usepackage{tikz-cd}
\usetikzlibrary{math}
\usetikzlibrary{calc}
\usepackage{stmaryrd} 
\usepackage{marginnote}
\usepackage{ifthen}
\usepackage[skip=10pt plus1pt, indent=10pt]{parskip}
\usepackage{cuted}

\usepackage{hyperref}

\newcommand{\N}{\mathbb{N}}

\newcommand{\R}{\mathbb{R}}

\DeclareMathOperator{\dist}{dist} 
\DeclareMathOperator{\trace}{trace} 
\DeclareMathOperator{\vol}{vol} 
\DeclareMathOperator{\dive}{div} 
\newcommand{\argument}{\mathord{\,\cdot\,}} 
\newcommand{\dx}{\;\mathrm{d}} 

\newcommand{\qad}{q_{\mathrm{ad}}}
\newcommand{\alphad}{\alpha_{\mathrm{ad}}}

\theoremstyle{definition}
\newtheorem{definition}{Definition}[section]
\newtheorem{remark}[definition]{Remark}

\newtheorem{assumption}[definition]{Assumption}

\theoremstyle{plain}
\newtheorem{proposition}[definition]{Proposition}
\newtheorem{lemma}[definition]{Lemma}
\newtheorem{theorem}[definition]{Theorem}
\newtheorem{corollary}[definition]{Corollary}

\numberwithin{equation}{section}

 \title{Dynamic access pricing control for fair and stable \\ resource sharing}
 \date{July 2026}

\author{C. King\protect\footnote{Department of Mathematics, Northeastern University, Boston, MA, USA}, F. R. Wirth\protect\footnote{Faculty of Computer Science and Mathematics, University of Passau, Passau, Germany}, H. Hamedmoghadam\protect\footnote{Dyson School of Design Engineering, Imperial College London, London, UK}, C. G. Cassandras \protect\footnote{Division of Systems Engineering, Boston University, Brookline, MA, USA}, \\ and R. Shorten\protect\footnote{Dyson School of Design Engineering, Imperial College London, London, UK}}

\begin{document}

\maketitle{}
          
\begin{abstract}
We consider the co-existence of users with different price sensitivities under dynamic (or surge) pricing that regulates how they access a shared resource.  To balance supply and demand in such settings, some dynamic pricing systems adopt a monotonically increasing price scheme which is well known to be socially regressive, driving out price-sensitive users in favor of price-insensitive ones.  Here, we propose a new dynamic pricing scheme based on a non-monotonic price function that allows all user classes to coexist more fairly when price-insensitive users are present, but automatically reverts to conventional dynamic pricing when they leave the system, even though the exact price-sensitivity composition is unobservable.  The resulting system consists of nonlinear ODEs, which we further investigate to characterize its equilibria and their domains of attraction, first for a single price-sensitive class competing with a price-insensitive one and then for the general $n$-class case. Numerical simulations demonstrate the effectiveness of the proposed design.
\end{abstract}


\section{Introduction}
Dynamic (or surge) pricing is a widely used strategy to manage user demand for a scarce shared resource/service during periods of congestion. Examples where pricing strategies have been deployed include highways, where price is sometimes used to maintain desirable flow rates, and communication networks where queue management strategies such as Random Early Detection (RED) are used to control transmission delays by maintaining small average queue lengths in congested routers \cite{bobfab,lehe2019downtown} (here price manifests as packets or ECN bits). While surge pricing has been used in many industries for some time, its popularity is now on the increase. For example, in a recent study in the Financial Times it was reported that up to one third of UK businesses plan to use dynamic pricing algorithms~\cite{ft_boe} and that this number is even higher in other countries~\cite{ft_surge}.

The historic literature on dynamic pricing is extensive. Much of the existing theoretical work concerns how a price should be set to maximize revenue~\cite{lin2006dynamic,den2015dynamic,chen2024deep,phillips2021pricing}. Other complementary work focuses 
on applications of dynamic pricing. For example, in the context of transportation, the pricing of vehicles accessing an urban area (also known as cordon pricing) dates back to 1975 in Singapore \cite{goh2002congestion} with the objective of managing congestion. Road and cordon pricing has also more recently been applied with the objective of reducing transport emissions \cite{green2020did,lehe2019downtown}. Further examples of dynamic pricing can be found in parking systems \cite{nourinejad2017impact} to reduce  demand for limited parking spaces, and also in energy applications~\cite{harsha2014optimal} where pricing is used as part of many demand-side management strategies (demand response). Dynamic pricing is also used in networking applications where a generalized notion of price is used to regulate access to internet routers \cite{jiang2008time,gizelis2010survey}, data pricing in machine learning pipelines \cite{cong2022data}, and very recently in the context of distributed ledgers (e.g., blockchain), where pricing acts as a mechanism to guarantee transactions are processed with low latency \cite{BobAnd}. In online platforms such as Uber ride-sharing and Airbnb lodging \cite{ubernature}, dynamic pricing is also used in the context of two-sided market design to both regulate access for a resource and incentivize additional supply. In all of these applications the fundamental assumption is that users are {\bf responsive} to price so that a supply-demand equilibrium emerges in the system under price control.

An emerging concern is that dynamic pricing may give rise to unfairness (and inequality)~\cite{eliasson2016congestion}. To understand this, it is worth noting that not all users acting under the influence of a surge pricing access control mechanism have the same level of price sensitivity. An extreme case in point is where one class of traffic (price-insensitive) does not respond to price at all, and another class (price-sensitive) does. Simply put, in such a \emph{heterogeneous} user population, price-sensitive users are disadvantaged: as the price increases, they will leave the system, making room for price-insensitive users. Policies of this nature could deny, for example, price-sensitive traffic access to services in key city zones, thereby contributing to \enquote{access poverty}, a term that is used to refer to the inability to access essential items and services. Moreover, when a large number of price-insensitive users compete for access, dynamic pricing may not even regulate demand, thereby negating the entire rationale for its use in the first place. In effect, the pricing strategy ends up merely \enquote{making space for the rich}.

One motivation for the use of dynamic pricing is to optimize revenue~\cite{phillips2021pricing}. However, while this applies to private resource providers, it may not apply to public institutions, such as cities and municipalities, whose objective is to provide access to a public resource (e.g., roadways, public parking facilities) with a high quality of service to users who support the resource through a pricing scheme. In such situations, ensuring fair and equitable access to a resource sometimes becomes an important factor for dynamic pricing algorithms, and standard dynamic pricing is not necessarily suited to achieving this goal.
Fairness concerns can also arise in the private sector. A recorded example of pricing-induced inequality in private resource-sharing is described in the Financial Times where the \enquote{ride-sharing app Uber refunded users in central London after its pricing engine briefly surged fares in the aftermath of the London Bridge terror attack in June 2017}~\cite{ft_surge}. Clearly in this  example, fair and equitable access to the resource was (in hindsight) deemed more important than maximizing revenue. This illustrates that even private service providers, thought to be purely profit-oriented, may need to abandon traditional surge pricing in some critical scenarios, even though that would lead to more revenue.

Our goal in this paper is to develop dynamic pricing algorithms that support the coexistence of heterogeneous user classes. A central challenge in this design is that the presence of a price-insensitive population is not directly observable, and may only be inferred from the system's aggregate response. This challenge is most visible in the case of \textit{bursty} inelastic demand, where price-insensitive users join and leave the system in response to exogenous events.
For example, in ride-hailing systems during periods of heavy rain, surges in price are driven by price-insensitive users who just want a cab to avoid becoming wet, and the pricing mechanism has no direct way to identify when this burst has subsided.
An intuitive strategy, aimed at making competition between the two classes fair, is to switch off dynamic pricing and set the price to zero.  Such approaches have been suggested in the past in TCP-based congestion control~\cite{bobrade1} in an attempt to enable coexistence of loss-based and delay-based TCP variants. However, dropping the price attracts a large volume of price-sensitive users, obscuring exactly the signal needed to detect the end of the burst, which makes it difficult to orchestrate switching between the pricing modes. A major contribution of this paper is a pricing mechanism that bypasses this fundamental \textit{observability} issue: when an inelastic burst is present, the proposed price function drops, allowing users with different price sensitivities to coexist; when the burst departs, the system automatically resumes surge pricing without the need to detect this transition explicitly. The mechanism builds on our prior work~\cite{bobrade1}, and involves designing queuing systems dynamics that resolve the observability issue described above, complementing related works that consider the impact of heterogeneous traffic classes in networks where some classes cannot be incentivized or are selfish \cite{yue,koll}.

The paper is organized as follows. First, in Section \ref{sec:predisc}, we establish a dynamic pricing system model which considers $n$ heterogeneous user classes requesting access to a shared resource. The classes differ in their price sensitivities, including a class which is entirely price-insensitive. We propose a price function which is \emph{non-monotonic} in the resource queue content, along with class-dependent sensitivity functions that model the dropout rates of different user classes in response to the price. In Section \ref{sec:normalmode}, we analyze the stability properties of this dynamic pricing system with a single price-sensitive class, establishing the existence of two equilibria (one stable and one unstable) and characterizing their domains of attraction. In Section~\ref{sec:multipleusers}, we extend the analysis to the general $n$-class case. Importantly, we show that the resulting dynamics allow heterogeneous user classes to coexist during periods when inelastic traffic is present, and that the system automatically returns to a surge pricing regime when the inelastic traffic leaves the system. In Section \ref{sec:interpretation} we then give an interpretation of the mathematical results we derive, namely that the impact of a surge of unresponsive users is transient. Section \ref{nozeroprice} then revisits the proposed pricing function and introduces a positive saturation level preventing it from reaching zero. 
Section \ref{sec:sim} presents simulation examples to illustrate the efficacy of our proposed solutions.

\section{Dynamic Pricing System Model}
\label{sec:predisc}
\begin{figure}[!ht]
    \centering
    \includegraphics[width=0.9\linewidth]{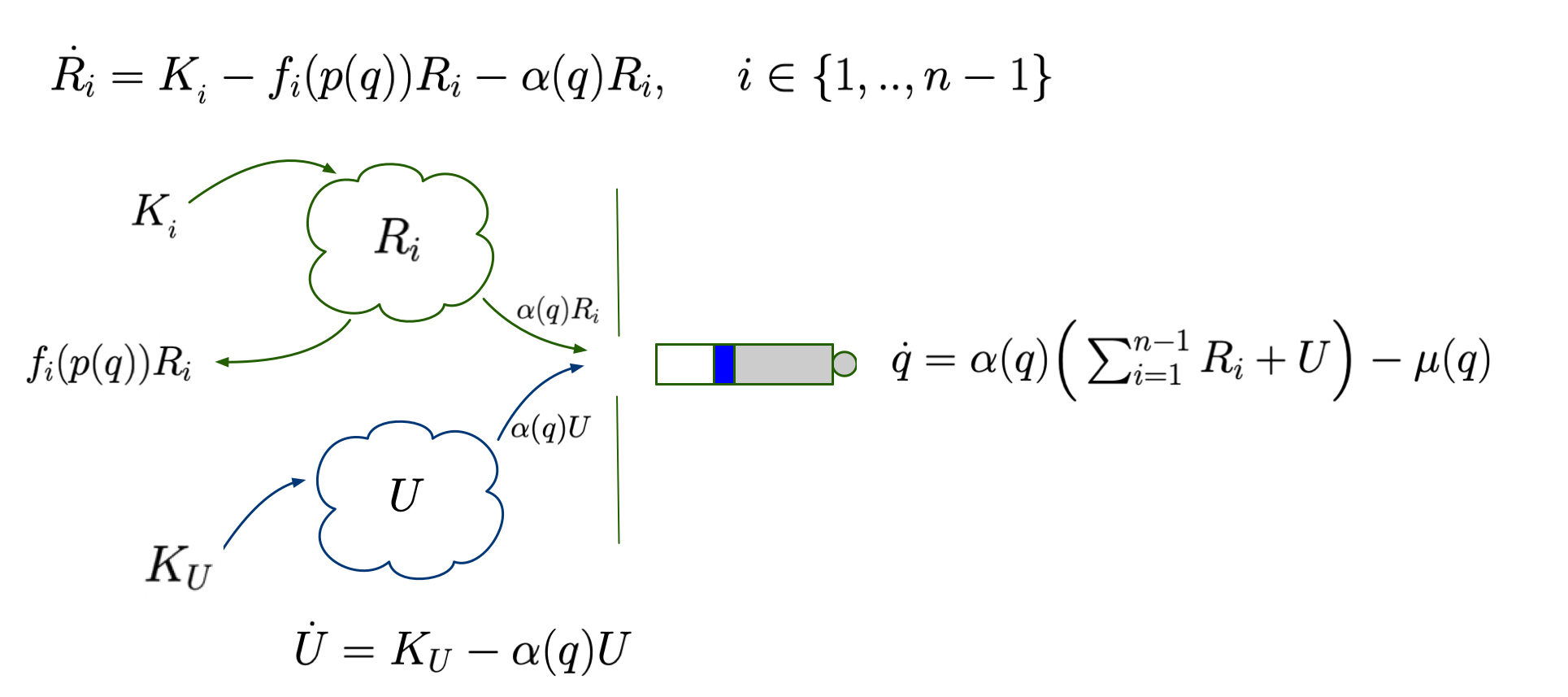}
    \caption{Price-sensitive and price-insensitive users competing for access to a resource under the influence of a price signal.
    }
    \label{fig:setting}
\end{figure}

Figure \ref{fig:setting} depicts a general setting in which $n-1$ ($n\geq 2$) price-sensitive user classes compete for access to a given resource, along with a price-insensitive one denoted by $U$. In the figure the $n-1$ price-sensitive classes are aggregated into one group, albeit with differing sensitivities. These classes are said to be {\em elastic} or {\em responsive} in the sense that their members are price-sensitive, and their respective queue contents are denoted by $R_i(t) \in \R_+$, with $i=1,\ldots,n-1$. On the other hand, the price-insensitive user class competes with the $R_i$ classes, and is \emph{inelastic} or {\em unresponsive}, in the sense that its members do not respond to the price signal. Its queue content is denoted by $U(t) \in \R_+$.
New users arrive into the $R_i$ and $U$ queues with rates $K_i(t)$ and $K_U(t)$, respectively, and leave the service queue at a rate given by  $\mu(q(t))$. In practice, $K_U(t)$ may be positive only over limited time intervals during which price-insensitive users compete with price-sensitive ones and achieving fairness becomes meaningful. The stability analysis that follows treats arrival rates as constant over sufficiently long intervals.

In this model, responsive users can leave their waiting queue in two ways. 
\begin{itemize}
    \item First, if the price is too high (in the case of price-sensitive users). This is governed by price sensitivity functions, defined in the sequel. 
    \item Second, if they accept the price and are admitted to the service queue. This is governed by an admission rate function, also defined in the sequel.
\end{itemize}
Letting $q(t) \in \R_+$ denote the length of the service queue at time $t$, we define a \emph{price function} $p(q(t))$ ($p:\R_+\to\R_+$, assumed to be locally Lipschitz continuous) which describes how the price charged for service evolves as a function of $q(t)$. The response of an individual class population to the price is defined by a \emph{price sensitivity function} $f_i(p(q(t)))$ ($f_i:\R_+\to\R_+$, $1\leq i \leq n-1$, assumed to be  continuously differentiable and non-decreasing).  The function $p(q(t))$ is always a design choice by the resource owner, whereas $f_i(p(q(t)))$ is normally not.  In many applications, this is an implicit property of a user population.

Entry to the service queue is governed by  an \emph{admission rate function} $\alpha(q(t))$ ($\alpha:\R_+\to\R_+$, assumed to be continuously differentiable and non-increasing). In many applications, this function characterizes the quality-of-service offered by the resource provider (which is in turn a function of the service queue length). In this paper we assume that dropouts are a function of price only, while admission to the service queue is controlled by $\alpha(q(t))$.

The flow dynamics of our system are given by (from now on the explicit dependency on  $t$ is dropped for notational ease): 
\begin{subequations}
\label{eq:system_bob}
\begin{align}
\label{eq:systemR12_bob}
    \dot R_i &= K_i - f_i(p(q)) R_i - \alpha(q) R_i,  \\
\label{eq:systemU1_bob}
    \dot U &= K_U - \alpha(q) U, \\
\label{eq:systemq1_bob}
\dot q &= \alpha(q) \left(\sum_{i=1}^{n-1}R_i + U\right)- \mu(q),
\end{align}
\end{subequations}
for $i=1,\ldots,n-1$. As seen in \eqref{eq:systemR12_bob}, the population of responsive users is reduced at a rate $f_i(p(q))R_i$, representing a proportion of users from class $i$ who choose to abandon access to the resource because of a high price. The population of price-responsive users is also reduced by the admission fraction $\alpha(q)R$ that actually accesses the service. The same reasoning applies to \eqref{eq:systemU1_bob} in which there are no price-dependent dropouts. Finally, in \eqref{eq:systemq1_bob} the dynamics of the resource occupancy depend on a service rate $\mu(q)$ (defined in the sequel) specifying the user departure process from this system. We write the above equations in a more compact form by 
extending the user class set to 
$i=1,\ldots,n$ and reserving the index $n$ to the $U$ class,  
\begin{subequations}
\label{eq:systemRq_bob}
\begin{align}
\label{eq:systemR13_bob}
    \dot R_i &= K_i - f_i(p(q)) R_i - \alpha(q) R_i,\\
\label{eq:systemq12_bob}
    \dot q &= \alpha(q) \sum_{i=1}^{n}R_i - \mu(q),
\end{align}
\end{subequations}
i.e., we set $f_{n}(p(q))=0$ and $R_{n} = U$. Together Eqs.~\eqref{eq:systemR12_bob},~\eqref{eq:systemU1_bob},~and~\eqref{eq:systemq1_bob} define a dynamic system of responsive, unresponsive, and service queues, noting, as mentioned earlier, that the presence of unresponsive traffic for certain applications is bursty and intermittent.

\begin{remark} 
\textbf{Dropout and wait mechanisms}
Observe that the model consisting of (\ref{eq:systemR12_bob}), (\ref{eq:systemU1_bob}), and (\ref{eq:systemq1_bob}) allows users to remain interested in accessing the resource until they are either admitted or (in the case of responsive users) they drop out. This model is rich enough to capture behavior found across many applications.
In this paper we limit ourselves to dropouts resulting from price sensitivity only, but
in more general settings, users may also decide to drop out as a function of service delay, or other factors.
\end{remark}

\begin{remark}\textbf{Role of pricing control.} In many applications, the purpose of pricing control is to protect the service to price-sensitive users by managing access to the service queue. This is only possible in the case when when unresponsive users are not competing with responsive ones. However, it is also true that in many situations when dynamic pricing is deployed, responsive and unresponsive traffic 
must coexist over short periods of time. For example, during bad weather, unresponsive traffic may enter an on-demand ride-hailing system. In such situations, the presence of unresponsive traffic may be viewed as a disturbance, and a principal contribution of this paper is to present \emph{mechanisms to ensure that the system returns to its natural equilibrium when the unresponsive traffic leaves the system}, even when the departure of this user class cannot be observed.
\end{remark}

In the remainder of this section, we discuss the pricing function $p(q)$, the price sensitivity functions $f_i(p(q))$, $i=1,\ldots,N$, the service rate function $\mu(q)$, and the admission rate function $\alpha(q)$. We adopt the point of view that the pricing function is (parametrically) specified and wish to study the equilibria (fixed points) in terms of the possible values of $q(t)$ that can be attained in this dynamic system. This is in contrast to an alternative viewpoint (the subject of ongoing research) in which first a desirable equilibrium $q^*$ is specified and we seek to determine a pricing policy that optimizes a given objective (e.g., deviation from $q^*$, or price) subject to system constraints (e.g., the service rate capacity of the resource). See \cite{christos} for an example of work in this direction.\newline

\subsection{The pricing function $p(q)$} 
Roughly speaking, two approaches are possible to adjust price to demand in a resource allocation problem. 
\begin{itemize}
    \item  One may adjust price in a PI control-like fashion, based on the difference between observed supply and demand \cite{schlote}. 
    \item A second approach is to adjust the price based on some quality of service metric, such as wait time or queuing delay. For example, Uber speaks of multipliers based on demand \cite{castillopaper}, which is one enunciation of this approach; in internet congestion control, Random Early Detection (RED) is another \cite{bobrade1}. This second approach only requires measurement of a queuing delay (or a proxy) and an appropriately specified pricing function for implementation. 
\end{itemize}
    We follow this latter approach in which a price signal $p(q)$ is used to modulate each user population $R_i$. Specifically, as mentioned, we assume members of $R_i$ leave their queue in a linear fashion as $f_i(p(q))R_i$. In RED queue management \cite{bobrade1}, for example, $p(q)$ is a non-increasing function of the length $q$. In what follows, we shall exploit a degree of freedom that exists by relaxing this assumption in order to realize queuing systems that have desirable dynamic properties. Our rationale for doing this, as opposed to the PI control-like approach, is as follows: 

    \begin{itemize}
        \item In many applications the purpose of the price is to deliver a good quality of service to users competing for a resource. Typically, this translates to a small service queuing delay (or small average service queue lengths). However, in many applications, the social cost of doing this is high, as price-sensitive traffic usually leaves the system to space way for traffic that is price-insensitive (i.e., financially disadvantaged users leave the system, hence relinquishing service to high-income or rich users). 
        \item If a situation prevails that a price function does not yield low average service queue lengths, or if there is a large volume of unresponsive traffic, then one may just as well set the price to zero and allow all traffic classes to compete for service in a fair manner (of course, at a cost of a low quality of service delivered to everyone). However, setting the price to zero beyond some threshold is also problematic, because all queues will 
        overflow
        making it impossible to observe from the queue dynamics when the nonresponsive user population leaves the system. A consequence of this is that developing a trigger to reactivate the dynamic pricing algorithm is difficult.
    \end{itemize}
    A price function of this kind may be defined as follows. Let $q_{m}$ 
    and $\beta$
    be positive constants and define, for $q\geq 0$. 
\begin{equation}
\label{eq:costfunction}
    p(q) = \left\{ \,\begin{matrix}
    \beta q & ,&  0\leq q \leq q_m, \\
    \beta ( 2 q_m -q) & ,& q_m\leq q \leq 2q_m, \\
     0 &,& 2q_m < q.
    \end{matrix} \right.
\end{equation}     
The key property expressed in \eqref{eq:costfunction} is that the price first increases monotonically, and then, beyond some threshold, decreases in a monotonic fashion. This simple property gives rise to unstable equilibria. As we shall see, these turn out to be beneficial in solving the coexistence problem of price-sensitive and insensitive users
and are an essential part of our design. For convenience, we have expressed these properties in a piecewise linear manner in \eqref{eq:costfunction}, but all arguments extend to more general price functions with these qualitative features. Note that we do not address how to select $q_{m}$.  This is the subject of extensions of this work, where we seek to optimize the pricing functions $f_i(\cdot)$ \cite{christos}.\newline 
\begin{remark} 
{\bf Revenue maximization and equity.} It could be argued that a pricing mechanism that involves reducing price is counter to the revenue-maximizing goal of a private operator. However, there are several reasons why incorporating fairness guardrails into the design of dynamic pricing algorithms. First, it is important to distinguish between private resource providers with such an optimization objective and ``public institutions'' (e.g., cities, municipalities) whose objective
is to provide access to a public resource (e.g., roadways, public parking facilities) with a high quality of service
level to users who support the resource through a pricing scheme. In the latter case, it is clear that fairness is
a key feature of ``high resource quality of service level''. As for the former case, it has recently been studied in \cite{christos}, where it is shown that a revenue-maximizing provider can still achieve this
objective during periods of time when a \emph{fairness constraint} is not violated. If the fairness constraint threshold
is set at a very low level, then the approach in \cite{christos} reduces to strict revenue maximization by combining both pricing
and access control. A second motivation for new types of dynamic pricing algorithms is that when  
price-insensitive users compete for access to a constrained resource, dynamic pricing alone will not be able to manage demand and some other access control mechanism is necessary \cite{policy_doc}. Third, regulatory concerns are now emerging around the impacts of surge pricing \cite{ft2,ft3}. For example, the UK watchdog launched a probe into the use of dynamic pricing for rock band Oasis's concerts that saw ticket prices skyrocket. This probe was centered on fairness issues that saw many fans priced out of tickets due to the presence of inelastic traffic.  Another regulatory concern is that the algorithms driving dynamic pricing can encourage implicit collusion between firms, raising prices overall and that dynamic pricing can also conceal price gouging in markets where there is fixed supply and little transparency\cite{ft3}. Finally, dynamic pricing could also be viewed as a primitive form of personalized pricing, which is seen as unfair in many jurisdictions and has been banned in some places (for example, Maryland)\cite{ft3}.
\end{remark}

\subsection{The class sensitivity functions $f_i(p(q))$}
Each class sensitivity function, $i \in \{1,...,n-1\}$, is assumed to be a continuously differentiable and nondecreasing function of $p(q)$. Furthermore, we allow 
the $n$'th function (corresponding to the unresponsive traffic) to be zero.
An example of such a function, $i \in 1,...,n-1$, is the linear price sensitivity function:
\begin{equation}\label{eqn:dropout_rate}
    f_i(p) = r_{i,1} p + r_{i,2},
\end{equation}
where \(r_{i,1}\) characterizes the price sensitivity of class \(i\), while \(r_{i,2} \ge 0\) captures its baseline tendency to leave the system regardless of price if \(r_{i,2} > 0\).

\subsection{The service rate function $\mu(q)$} 
For simplicity, here, we assume a service rate
$\mu(q)$ that increases linearly in $q$ and levels off at a certain stage, i.e., there are constants $\mu^*,q_c>0$ such that $\mu:\R_+\to \R_+$ is given by
\begin{equation}
    \label{eq:mudef}
    \mu(q) = \left\{ \begin{matrix}
       \displaystyle{\frac{\mu^*}{q_c} q} && 0\leq q \leq q_c, \\[1ex]
     \mu^* && q_c \leq q.\newline
    \end{matrix}
    \right.
\end{equation}

\subsection{The admission rate function $\alpha(q)$} 
Aside from reactions to pricing,
the rate at which the population $R_i$ grows is also influenced by a term $\alpha(q)$ as seen in Eqs. 
(\ref{eq:systemR12_bob}) and (\ref{eq:systemU1_bob}). We use this term to model users who find the price acceptable and for whom the offered quality of service is satisfactory. Similar to the price function, we assume that users are admitted to the service queue at a rate $\alpha(q)R_i$, and that $\alpha(q)$ is a non-increasing function of $q$. We shall discuss $\alpha(q)$ in more detail when we address the equilibrium states of our model. 

\section{A single price-responsive class}
\label{sec:normalmode}
We begin the analysis of the dynamic pricing system model (\ref{eq:systemR12_bob}), (\ref{eq:systemU1_bob}), and (\ref{eq:systemq1_bob}) by considering the case of a single price-sensitive user class ($n=1$). This will allow us to extract the salient features of this system, in particular understanding the behavior of the system when it starts just after the insensitive class has left the system and and attains an equilibrium, if one exists. Subsequently, in Section \ref{sec:multipleusers} we will extend the analysis to the general case of $n$ user classes.

We are interested in the situation that prevails when a single elastic population competes with an inelastic population for a shared resource, just after the inelastic population has left the system at a time $t^*$. Our specific interest is to understand what happens starting from the initial state at $t^*$.
In this case, the dynamics reduce to, with the subscript $i$ removed for clarity:
\begin{subequations}
\label{eq:system}
\begin{align}
\label{eq:systemR}
    \dot R &= K_R - f(p(q)) R - \alpha(q) R,\\
\label{eq:systemq}
    \dot q &= \alpha(q) R - \mu(q)
\end{align}
\end{subequations}
with the explicit assumption that $K_R>0$ is a constant user inflow 
over a sufficiently long time interval to allow the system to reach an equilibrium, if one exists. 

Where convenient we will abbreviate the right hand side of $\eqref{eq:system}$ by $F$, leading to the equivalent formulation
\begin{equation}
    \dot {\begin{bmatrix}
        R \\ q
    \end{bmatrix}} = F(r,q).
\end{equation}
We will assume $K_R> \mu^*$, so that the arrival rate for the demand queue exceeds the service rate for the service queue. Otherwise, it is clear that the second queue will always be empty so there are no interesting dynamics and no real requirement for an access policy to the resource. 
The function $\alpha:\R_{\geq 0}\to \R_{>0}$ will be chosen as a continuously differentiable function, to determine the desired fixed points of our system. 
As we will see, the exact nature of $\alpha(q)$ critically affects the number of fixed points in our system. Note also that for any $\alpha(q)$ of this type, 
$(R^*,q^*) \in \R^2_{\geq 0}$ is a fixed point of \eqref{eq:system} if and only if
\begin{equation}
\label{eq:fp-conds3}
    R^* = \frac{\mu(q^*)}{\alpha(q^*)}, \quad \frac{K_R-\mu(q^*)}{\mu(q^*)} = \frac{f(q^*)}{\alpha(q^*)}.
\end{equation}

We aim to have two fixed points corresponding to the low and high congestion regimes of the queues, as defined through Eq. \eqref{eq:costfunction}. Our rationale for targeting two points is to use one of the fixed points to give a desired level of performance when the price-unresponsive (inelastic) traffic is not present, i.e., the stable equilibrium; and to use the other, i.e., unstable equilibrium, as a \enquote{spring} to push the system back to the desired stable equilibrium when the price-unresponsive traffic leaves the system. If the desired equilibrium prices of our system are $0<p_1< \beta q_m$ in the low congestion regime, in which $q \in [0,q_m)$ and $0<p_2<\beta q_m$ in the high congestion regime in which $q \in (q_m, 2 q_m]$, then we obtain for the corresponding fixed points $(R^s,q^s)$ and $(R^u,q^u)$  the conditions that
\begin{align}
\label{eq:fp1_condition}
\begin{split}
q^s &= \frac{f^{-1}(p_1)}{\beta}, \quad
\alpha(q^s) = f(p_1)\frac{ \mu(q^s)}{K_R-\mu(q^s)}, \\
R^s &= \frac{K_R-\mu(q^s)}{f(p_1)},
\end{split}
\end{align}
and
\begin{align}
\label{eq:fp2_condition}
\begin{split}
q^u &= 2 q_m - \frac{f^{-1}(p_2)}{\beta},\quad
\alpha(q^u) = f(p_2) \frac{ \mu(q^u)}{K_R-\mu(q^u)},\\
    R^u &= \frac{K_R-\mu(q^u)}{f(p_2)}.
\end{split}
\end{align}

It will be convenient to consider the nullclines, i.e. the curves along which $\dot R = 0$ and $\dot q = 0$. These are given by the points $(q,\eta_i(q))$ with
\begin{align}
    \label{eq:eta1def}
    \dot q = 0 \quad &\Leftrightarrow \quad R(q) = \eta_1(q) := \frac{\mu(q)}{\alpha(q)},\\
    \dot R = 0 \quad &\Leftrightarrow \quad R(q) = \eta_2(q):= \frac{K_R}{\alpha(q)+ f(p(q))}.
    \label{eq:eta2def}
\end{align}
With this, fixed points are precisely the zeros of the function
\begin{equation}
    g(q) := \eta_2(q) - \eta_1(q) = \frac{K_R}{\alpha(q)+ f(p(q))} - \frac{\mu(q)}{\alpha(q)}.
\end{equation}

To avoid the existence of further fixed points, all that is required is that $\alpha(q)$ is chosen such that the second condition in \eqref{eq:fp-conds3} is not satisfied for $q\neq q^s, q^u$.  We formulate this as an explicit condition on $\alpha(q)$. 
\begin{assumption}
    \label{ass:alphaprops}
    Consider system \eqref{eq:system} with a cost function given by \eqref{eq:costfunction}.
    Let the sensitivity function $f:\R_+\to \R_+$ be continuously differentiable with strictly positive derivative.
    Consider a fixed maximal queue length $q_{\mathrm{max}} \geq 2 q_m$.
    We call a continuously differentiable admission rate $\alpha: \R_+\to [0,\infty)$ \emph{admissible}, if 
    \begin{enumerate}[(i)]
        \item $\alpha(q)$ is positive and strictly decreasing on $[0,q_{\mathrm{max}})$.
        \item $\alpha(q) = 0$, if $q\geq q_{\mathrm{max}}$.
        \item there are exactly two points $q^s,q^u\in (0,\infty)$ solving the equation
    \begin{equation}
     \label{eq:fixed-point-conds}
\eta_1(q) = \eta_2(q),
    \end{equation}
    where $q^s \in (0,q_m)$, $q^u\in (q_m, 2q_m)$ and $\eta_1'(q^u) < \eta_2'(q^u)$. Furthermore, $\mu(q^s) = \mu(q^u) = \mu^*$.
    \end{enumerate}
\end{assumption}

By \eqref{eq:eta1def}, \eqref{eq:eta2def}, condition \eqref{eq:fixed-point-conds} defines the fixed points of system~\eqref{eq:system}. In our model the situation of interest is that the service queue is operating at full capacity at the fixed points. We have made this explicit in the final point of Assumption~\ref{ass:alphaprops}.

For the analysis of the dynamic behavior we first recall the concept of forward (or positive) invariance. Consider a differential equation $\dot x=h(x)$ defined on $D\subset \R^n$ with unique solutions $\varphi(\cdot,x_0)$ to initial value problems $x(0)=x_0$ existing on $\R_+$. A subset $C\subset D$ is called \emph{(forward) invariant}, see e.g. \cite{yorke1967invariance}, or \cite[Theorem 4.3.8]{ClarLedy98}, if for all $x\in C$ the corresponding solution $\varphi(t,x)\in C$ for all $t\geq0$.
A sufficient condition for forward invariance of $C$ is that the vector field $h$ is subtangential at the boundary of $C$, \cite[Theorem 3.3]{yorke1967invariance}. In case $C$ is a polytope, it is sufficient that for all facets $K$ of $C$, i.e. $n-1$ dimensional intersections of the boundary with a hyperplane, it holds for the outward normal $n$ of $K$, that $\langle n, h(x) \rangle \leq 0$ for all $x\in K$. This is the condition that we will frequently use in the sequel.
A basic modeling requirement is that the queue lengths can never become negative. 
It is thus of interest to note that our state space is {\bf forward invariant} under the dynamics of system \eqref{eq:system} as shown next.

\begin{proposition}
    \label{prop:invariance}Under Assumption~\ref{ass:alphaprops} and
    for system \eqref{eq:system} we have
    \begin{enumerate}[(i)]
        \item the positive orthant $\R^2_+:=\{ x \in \R^2 \;;\; x_i \geq 0, i=1,2\}$ is forward invariant.
        \item the set $[0,\infty) \times [0, q_{\mathrm{max}}]$ is forward invariant.
    \end{enumerate}
\end{proposition}

\begin{proof}
    (i) If $x=(0,q)$, $q\geq 0$ we have $\dot R = K_R> 0$
    from (\ref{eq:systemR}). If $x=(R,0)$, $R\geq0$, we have $\dot q = \alpha(0)R\geq0$
    from (\ref{eq:systemq}). This implies that for all points $x$ on the boundary of $\R^2_+$ (with the exception of $0$) the inner product of outward normal in $x$ is negative; i.e., the vector $\dot{x}$ points into the region $\R^2_+$.
    This implies invariance.

    (ii) By (i) we only need to check the behavior of the flow on the boundary section $\{(R,q_{\mathrm{max}}); R \geq 0\}$. In these points we have by Assumption~\ref{ass:alphaprops}\,(ii) that $\dot q = -\mu(q_{\mathrm{max}})<0$. The assertion follows as in (i).
\end{proof}We then have the following lemma on the {\bf local stability analysis}.

\begin{lemma}
\label{lemma:localstability-2d}
    Consider the system \eqref{eq:system} and let Assumption \ref{ass:alphaprops} hold.
    The right hand side of Eqs. \eqref{eq:system} are continuously differentiable in the fixed points $x^s=(R^s,q^s)$ and $x^u=(R^u,q^u)$. In every point $x$ of differentiability of $F$ the Jacobian\footnote{We have omitted the arguments $x=(R,q)$ in the presentation of the Jacobian to avoid overloaded notation. Also note that $(f\circ p)':=\dx (f\circ p)/\dx q$, and $\alpha':=\dx \alpha/\dx q$, etc.} is
 \begin{equation}
 \label{eq:2dim-linearization}
     DF(x) = \begin{bmatrix}
         - (f\circ p+\alpha) & - R((f\circ p)' + \alpha')\\
         \alpha & R\alpha' - \mu'
     \end{bmatrix}.
 \end{equation} 
    The divergence of the vector field in its points of differentiability is
    \begin{align}
    \label{eq:divergence}
        \mathrm{div} F(R,q) = 
        \left\{ 
        \begin{matrix}
            - f(p(q)) - \alpha(q) + \alpha'(q)R - \frac{\mu^*}{q_c} &
             & 0< q < q_c, \\
            - f(p(q)) - \alpha(q) + \alpha'(q)R
             & & q_c < q.
        \end{matrix}
        \right. 
    \end{align} 
    In particular, it follows that
    \begin{enumerate}[(i)]
    \item  the fixed point $(R^s,q^s)$ is locally asymptotically stable;
    \item 
    the fixed point $(R^u,q^u)$ is unstable with a linearization with one positive and one negative eigenvalue;
    \item The system \eqref{eq:system} does not have nontrivial periodic solutions in the interior of $\R_{\geq0}^2$.
\end{enumerate}
\end{lemma}

\begin{proof} 
For computational details concerning the Jacobian and the divergence we refer to the Appendix~\ref{appendix}. Note that \eqref{eq:2dim-linearization} follows directly from \eqref{eq:system}. For \eqref{eq:divergence}, recall that the divergence is the trace of the Jacobian and use \eqref{eq:mudef}.
    Statements (i) and (ii) are immediate consequences of the more general result presented in Theorems~\ref{thm:musers-localstability} and~\ref{thm:musers-localinstab}. We include proofs here for completeness. Note that the assumption that $\alpha(q)$ is nonincreas\-ing implies that $\alpha'(q) \leq 0$ for all $q\geq 0$. 
 
    (i) The trace/determinant criterion for the Hurwitz property of $A\in \R^{2\times 2}$ states that $A$ is Hurwitz if and only if $\trace (A) <0$ and $\det(A) >0$. For the matrix $A$ in question and $x^s = (R^s,q^s)$ we have 
    \begin{align}
        \trace& (DF(x^s))= - (f(p(q^s)) + \alpha(q^s)) + R^s\alpha'(q^s) - \mu'(q^s) < 0,\notag
    \end{align}
    where we have used $f(p(q^s)),\alpha(q^s),R^s>0$, $\mu'(q^s)\geq 0$ and $\alpha'(q^s)\leq 0$. Moreover,
    (again dropping the argument $q^s$ for simplicity)
       \begin{multline}  
\nonumber
     \det(DF(x^s)) = -(f\circ p+\alpha)(R^s\alpha'- \mu') + R^s((f\circ p)' + \alpha')
\alpha \\
       = - (f\circ p)R^s\alpha' + (f\circ p+\alpha)\mu' + R^s(f\circ p)' \alpha >0,
       \label{eq:detcalchelp}
       \end{multline}
    where we have used the assumption that $0\leq q^s< q_m$, so that $(f\circ p)'(q^s) >0$. As the linearization in the fixed point $x^s =(R^s,q^s)$ is Hurwitz it follows from Lyapunov's linearization theorem \cite{khalil2002nonlinear}
    that the fixed point $x^s$ is asymptotically stable for system \eqref{eq:system}.

    (ii) For $x^u = (R^u,q^u)$, we have as before
    \begin{align}
        \trace (DF(x^u)) = 
        - (f(q^u) + \alpha(q^u)) + R^u\alpha'(q^u) - \mu'(q^u) < 0,\notag
    \end{align}
    so that at least one of the eigenvalues of $DF(x^u)$ has negative real part. In addition,
    we have
        \begin{equation}
    \label{eq:2dim-unstabledet}
        \det(DF(x^u))  
        = -(f\circ p+\alpha)(R^u\alpha'- \mu') + R^u((f\circ p)' + \alpha')\alpha.
     \end{equation}
    From the condition $\eta_1'(q^u) < \eta_2'(q^u)$ of Assumption~\ref{ass:alphaprops}\,(iii) we have with \eqref{eq:eta1def}, \eqref{eq:eta2def} that
    \begin{equation}
    \label{eq:543}
        \frac{\mu'}{\alpha} - \eta_1(q^u) \frac{\alpha'}{\alpha}= \frac{\mu' \alpha - \mu \alpha'}{\alpha^2} = \eta_1' < \eta_2' 
        = -\frac{K_R((f\circ p)' + \alpha')}{(f\circ p)^2} = -\eta_2(q^u) \frac{(f\circ p)' + \alpha'}{f\circ p}.
        \end{equation}
    Multiplying the leftmost and rightmost expressions of \eqref{eq:543} with $\alpha(\alpha + f\circ p)$ and using $R^u= \eta_1(q^u) = \eta_2(q^u)$, we obtain
    \begin{equation}
       (\mu' - R^u \alpha') (\alpha+f\circ p) < - R^u((f\circ p)' + \alpha')\alpha.
    \end{equation}
    A direct comparison with \eqref{eq:2dim-unstabledet} shows that this implies $\det(DF(x^u)) <0$. Thus the linearization in $x^u$ has a real, positive eigenvalue and the fixed point is unstable for system \eqref{eq:system}.

    (iii) Assume that $\psi$ is a nontrivial periodic solution of \eqref{eq:system} lying in the interior of $\R^2_{\geq 0}$. Then the orbit $\{ \psi(t) \;;\; t\geq 0 \}$ is a Jordan curve that separates the bounded interior $U$ of the orbit from the exterior, \cite[Theorem~VII.4.1]{hartman1982ordinary}. For the flow $\varphi$ generated by \eqref{eq:system} it follows by uniqueness of solutions for all $t\geq 0$ that $\varphi_t(U) = U$, and so in particular the volume of $\varphi_t(U)$ is constant. On the other hand, from \eqref{eq:divergence} we have $\dive F(R,q) < 0$ for all $(R,q)$ with $R>0,q>0$ because $\alpha'(q) \leq 0$. This holds in particular for all $(R,q) \in U$. It follows from the divergence theorem for Lipschitz continuous vector fields, see \cite[Proposition 1]{bessa2017flowbox}, that $\vol (\varphi_t(U)) < \vol (U)$ for all $t>0$. This contradiction shows that a periodic solution $\psi$ does not exist.
\end{proof}

\subsection{The domain of attraction}
\label{sec:domain}
As we have seen from our local stability analysis, under Assumption~\ref{ass:alphaprops}, there is always an asymptotically stable fixed point $x^s=(R^s,q^s)$ with $0< q^s < q_m$. We now aim to provide some analytical insights into the domain of attraction of this fixed point. Recall that the domain of attraction is defined as the set of initial conditions $\mathcal{A}(x^s)$ from which the trajectory $\varphi(t;x)$ converges to the fixed point, i.e. we have
\begin{equation*}
    \mathcal{A}(x^s) = \{ x \in \R^2_+ \;;\; \lim_{t\to \infty} \varphi(t;x) = x^s \}.
\end{equation*}

\begin{figure}[t]
    \centering
    \includegraphics[width=0.7\linewidth]{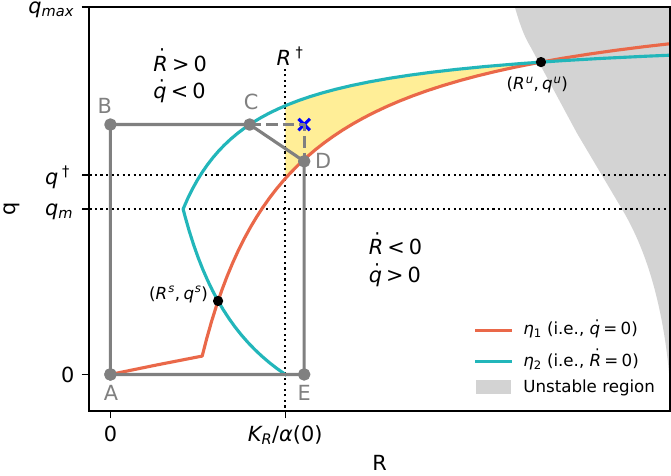}
    \caption{Sketch of the domain of attraction $\mathcal{A}$ of the stable fixed point $x^s=(R^s, q^s)$.
    Here arbitrary values are chosen for the system parameters, and an appropriate admission function of the form $\alpha(q)=\max(0, c_1q+c_2)$ is used with $c_1$ and $c_2$ calculated according to the conditions in Eqs.~\eqref{eq:fp1_condition} and \eqref{eq:fp2_condition}. The hatched area in yellow shows the set of all points with a given $q \in \left(q^\dagger, q^u\right)$ 
    and respectively $R \in (\max\{R^\dagger,\eta_2(q)\}, \eta_1(q)]$, addressed in Theorem~\ref{theo:domain}. The area shaded in gray marks the unstable region.
    }
    \label{fig:domain}
\end{figure}

The geometric construction is represented in Fig.~\ref{fig:domain}. In the $(R,q)$ plane we draw the curves $\eta_1,\eta_2$ defined in \eqref{eq:eta1def}, \eqref{eq:eta2def} of points for which $\dot R = 0$ and $\dot q = 0$. 

To provide the analytic results for the domain of attraction we note the following simple facts.

\begin{itemize}
    \item The curves $\dot q = 0$ and $\dot R = 0$ are given as graphs of continuous functions of $q$, which we have denoted $\eta_1,\eta_2$. By our assumption these graphs intersect in exactly two points, namely the fixed points $x^s=(R^s,q^s)$ and $x^u=(R^u, q^u)$.
    \item the function $\eta_1$ is increasing as a function of $q$, as $\mu$ is increasing and $\alpha$ is decreasing, see \eqref{eq:eta1def}.
    \item The function $\eta_2$ is increasing on the interval $[q_m,q_{\max}]$ as $\alpha$ and $f$ are both decreasing there. 
    We define 
    \begin{align}
    \label{eq:Rmaxdef}
            R^{\dagger} := \max \{ \eta_2(q) \;;\; q \in [0,q_m] \} 
        \geq \max \left\{ \frac{K_R}{\alpha(0)}, R^s \right\}
        \end{align}
    and let $q^\dagger := \eta_1^{-1}(R^\dagger)$.
\item We have $\eta_1(0) = 0$ and $\eta_2(0) = \frac{K_R}{\alpha(0)}$. Also with the intersection points $(R^s,q^s),(R^u,q^u)$, it may be seen that
    $\eta_1(q) < \eta_2(q)$, $q\in [0,q^s)$, and $\eta_2(q) < \eta_1(q)$, 
    $q\in (q^s,q^u)$, and finally 
    $\eta_1(q) < \eta_2(q)$, 
    $q \in (q^u,q_{{\mathrm{max}}})$.
\end{itemize}

With these observations we can prove the following theorem.

\begin{theorem}
    \label{theo:domain}
If Assumption~\ref{ass:alphaprops} holds and $R^u > R^{\dagger}$, then for any $\hat{q}$ satisfying
\begin{equation}
    \hat{q} \in \left(q^\dagger, q^u\right)
\end{equation}
and any $\hat{R} \in (\max\{R^\dagger,\eta_2(\hat{q})\}, \eta_1(\hat{q})]$, the interior of the polygon defined as the convex hull of the points, $A = (0,0)$, $B = (0,\hat{q})$, $C = (\eta_2(\hat{q}),\hat{q})$, $D = (\hat{R},\eta_1^{-1}(\hat{R}))$, $E = (\hat{R},0)$,
is a forward invariant set for \eqref{eq:system} that is contained in $\mathcal{A}(x^s)$.
\end{theorem}

\begin{proof}
   It is sufficient to show that the polygon $P = \mathrm{conv}\{A,B,C,D,E\}$ is forward invariant. Indeed, the fixed point $x^s$ is the only fixed point in $P$ and by Lemma~\ref{lemma:localstability-2d} the system does not have nontrivial periodic solutions. 
   By the invariance of $P$, the $\omega$-limit set $\omega(x_0)$ of every initial condition $x_0\in P$ lies in $P$.
   By the Poincar{\'e}-Bendixson theorem, \cite[Theorem~VII.4.1]{hartman1982ordinary}, 
   $\omega(x_0)$ is a periodic solution or contains a fixed point. As $x^s$ is the unique fixed point in $P$ and asymptotically stable, we have $\lim_{t\to\infty} \varphi(t,x_0) = x^s$ for every $x_0 \in P$.

   To show invariance it is sufficient to consider the segments between the vertices of the polygon. 
   Note that from the choice of $\hat{q}> q^\dagger$ we have $\eta_1^{-1}(K_R/\alpha(0))< \hat{q} < q^u $. This implies with the choice of $\hat{R}$ that $ \max\{K_R/\alpha(0), \eta_2(\hat{q})\} < \hat{R} \leq \eta_1(\hat{q}) < \eta_1(q^u) = R^u$. 
   It is easy to see that on the segment $\overline{AB}:= \{\lambda A + (1-\lambda)B; \lambda \in (0,1)\}$ we have $\dot R =K_R > 0$; on $\overline{BC}$ it holds $\dot q < 0$ by definition of $\eta_1$; on $\overline{DE}$ it holds that $\dot R <0$ by definition of $\eta_2$ and as $\hat{R}> R^\dagger$; and on $\overline{EA}$ we have $\dot q = \alpha(0) R > 0$, unless $R=0$. 

   The only interesting segment is therefore the segment $\overline{CD}$. For this note that by construction $ \eta_2(\hat{q})< \hat{R}$ and $\hat{q} > \eta_1^{-1}(\hat{R})$. Thus the vector $C-D$ is strictly negative in the first and strictly positive in the second component. Consequently the outside normal $v$ to $P$ on the edge $\overline{CD}$ is a positive vector in both components. The segment $\overline{CD}$ lies entirely in the region in which $\dot R <0$ and $\dot q <0$ (the shaded area in Fig.~\ref{fig:domain}). Thus for any $x \in \overline{CD}$ we have
   \begin{equation*}
       \langle v , F(x) \rangle < 0.
   \end{equation*}
   This shows that trajectories of \eqref{eq:system} cannot leave $P$ by passing through $\overline{CD}$. By continuity of the flow it is not necessary to check the vertices of the polytope $P$, so the proof of invariance is complete.
\end{proof}

\begin{figure}[t]
    \centering
    \includegraphics[width=0.75\linewidth]{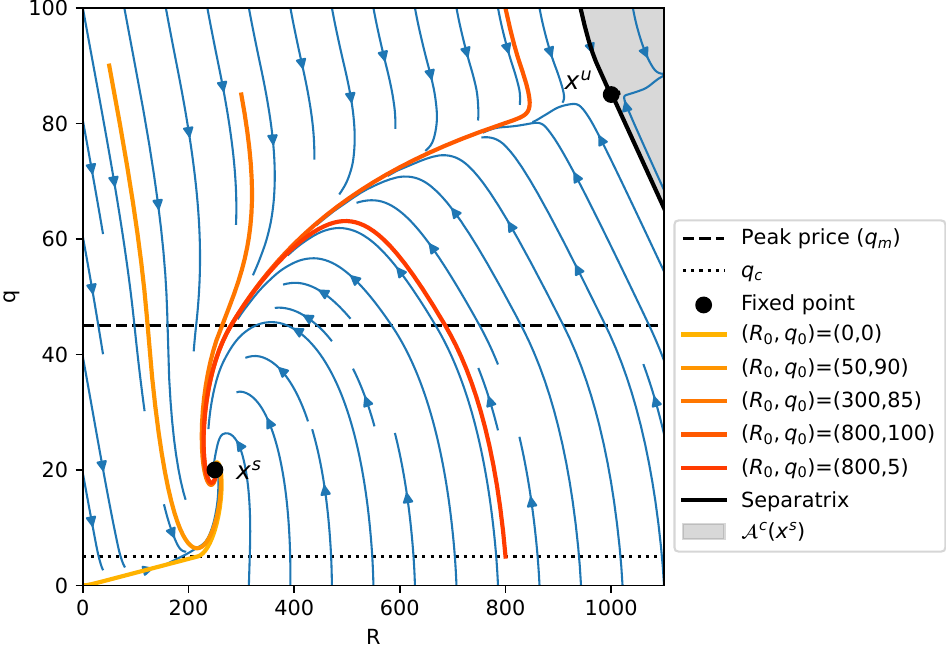}
    \caption{Phase portrait for a system based on the model in Eqs.~\eqref{eq:system}. The system fixed points, the inherent system constant $q_c$, the pricing design constant $q_m$, and trajectories for multiple initial states are plotted on top of the phase plane. The direction and color of arrows show the phase and magnitude at each point; with the magnitude gradually decreasing from red to yellow, green, and blue. }
    \label{fig:phaseplane2D}
\end{figure}

Figure~\ref{fig:phaseplane2D} depicts the phase portrait of the system defined in Eqs.~\eqref{eq:system}, marking the coordinates of the two fixed points and several trajectories for different initial conditions.

\subsection{Modified admission policy}
\label{subsec:modified-1user}

\begin{figure}[htp]
    \centering
    \includegraphics[width=0.6\linewidth]{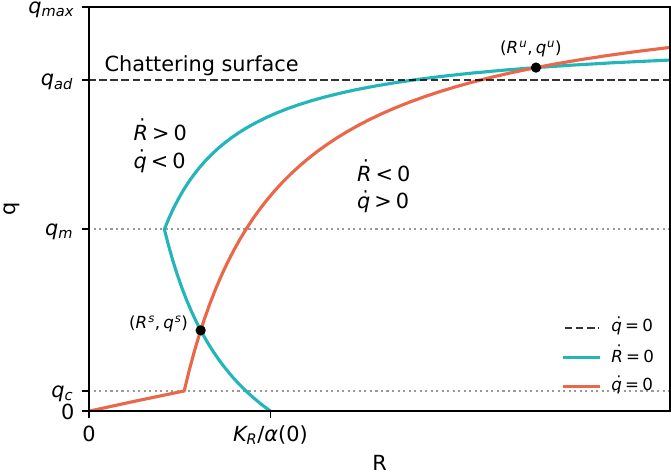}
    \caption{Sketch of the chattering surface in the $(R, q)$ space.
    }
    \label{fig:domain2}
\end{figure}

Figure \ref{fig:domain} depicts the domain of attraction for our queuing system. We note that within the gray region there exists a subregion where both $\dot{q}>0$ and $\dot{R} >0$. 
 Clearly this region is problematic as both state variables grow in it. To overcome this problem, and to yield a queuing system that is globally attractive, we now make a final modification to system dynamics. In particular, we define a maximum queue length to create a chattering boundary along which $\dot{R} < 0$ and $\dot{q}=0$ and along which the system will converge toward the domain of attraction of the stable equilibrium point as $\dot{R}$ is bounded away from zero along this surface. This is depicted in Fig.~\ref{fig:domain2}. 

To achieve the desired behavior, we introduce an admittance bound for the service queue, which we denote by $q_{\mathrm{ad}}$. At and beyond this bound the number of users that can enter the service queue per time period is bounded by the number of users that can be served in the same time period.  We  choose
$q_{\mathrm{ad}}$ to be in the interval $(q^\dagger, q^u)$, with $q^\dagger$ defined in \eqref{eq:Rmaxdef}.

The admission policy to the service queue is now changed from the original $\alpha(q)$ to 
\begin{equation}  \label{eq:chattering_queue_admittance}
     \alpha_{\mathrm{ad}}(q,R)= \left\{ \begin{matrix}
         \alpha(q) R & \quad & 0 \leq q < q_{\mathrm{ad}}\,,\\
         \min \{ \alpha(q) R, \mu^* \} && q \geq q_{\mathrm{ad}}
     \end{matrix}\right.
\end{equation}
With this new admission policy the dynamics formally change to 
\begin{subequations}
\label{eq:system-chatter}
\begin{align}
\label{eq:systemR-chatt}
    \dot R &= K_R - f(p(q)) R - \alpha_{\mathrm{ad}}(R,q),\\
\label{eq:systemq-chatt}
    \dot q &= \alpha_{\mathrm{ad}}(R,q) - \mu(q).
\end{align}
\end{subequations}
Note that in the region $[0,\infty) \times [0,q_{\mathrm{ad}})$ the dynamics in \eqref{eq:system-chatter} coincide with those defined in \eqref{eq:system}. On the chattering surface depicted by a dashed line in Fig.~\ref{fig:domain2} we have introduced a discontinuity by setting $\dot q = 0$ for all points $(R,q)$ with $q\geq q_{\mathrm{ad}}$ and $\alpha(q) R > \mu^*$.

As we are now dealing with a differential equation with discontinuous right hand side, some care is required concerning the solution concepts. We refer to \cite{Fili88} for details on this. Here we just mention that the two concepts of interest are Carath{\'e}odory solutions\footnote{Absolutely continuous functions, that satisfy the right hand side almost everywhere.} and Filippov solutions\footnote{Solutions of the differential inclusion obtained by Filippov regularization.}. It is not hard to see that the modification is sufficiently benign, so that for every initial condition $(R_0,q_0) \in [0,\infty) \times [0, q_{\mathrm{ad}}]$ we have unique Carath{\'e}odory solutions of \eqref{eq:system-chatter} in forward time that are defined on the time interval $[0,\infty)$ and the sets of Carath{\'e}odory and Filippov solutions coincide. 

In the situation of system \eqref{eq:system-chatter} together with \eqref{eq:chattering_queue_admittance} we have the following stability result, recalling the constants $R^{\dagger}$, $q^{\dagger}$ defined in \eqref{eq:Rmaxdef}.

\begin{boxedminipage}{0.98\textwidth} 
\begin{theorem}{\bf Main Result (single responsive class)}
\label{th:main_result_sec3}
Let Assumption~\ref{ass:alphaprops} hold.
    Consider \eqref{eq:system-chatter} with \eqref{eq:chattering_queue_admittance} and a continuously differentiable admission rate $\alpha(q)$.
Assume $R^u > R^{\dagger}$ and
$q^\dagger <q_{\mathrm{ad}} < q^u$.
    Then $x^s=(R^s,q^s)$ is a locally asymptotically stable fixed point of \eqref{eq:system-chatter} and for every initial condition $x=(R,q) \in X_{\mathrm{ad}}:=[0,\infty) \times [0,q_{\mathrm{ad}}]$ we have for the trajectory $\varphi(t;x)$
    \begin{equation*}
        \lim_{t\to \infty} \varphi(t;x) = x^s.
    \end{equation*}
\end{theorem}
\end{boxedminipage}

\begin{proof}
    The local stability of $x^s$ was already shown in Lemma~\ref{lemma:localstability-2d}. For the remainder of the proof we distinguish two cases: (i) initial conditions $(R_0,q_0) \in X_{\mathrm{ad}}$ such that the corresponding solution of \eqref{eq:system-chatter} does not intersect the chattering surface; (ii) all other initial conditions $(R_0,q_0) \in X_{\mathrm{ad}}$. Also observe that $X_{\mathrm{ad}}$ is trivially forward invariant under \eqref{eq:system-chatter}.
    
    (i) All trajectories $x(\argument)=(R(\argument),q(\argument))$ of \eqref{eq:system-chatter} starting in $X_{\mathrm{ad}}$ that do not intersect the chatter line coincide with trajectories of the system \eqref{eq:system}. For such trajectories $q(t) < q_{\mathrm{ad}}$ for all $t\geq 0$  and in this case $q(\argument)$ is bounded by assumption. Also $\dot R (t) <0$ for all $t$ where $R(t)$ is sufficiently large. Thus $(R(\argument),q(\argument))$ is a bounded trajectory and hence has a nonempty $\omega$-limit set. By Lemma~\ref{lemma:localstability-2d} and the Poincar{\'e}-Bendixson theorem, this has to be a fixed point, whence $\omega(x) =\{ x^s\}$ by forward invariance of $X_{\mathrm{ad}}$.

    (ii) If $(R_0,q_0)\in X_{\mathrm{ad}}$ is such that the corresponding trajectory of \eqref{eq:system-chatter} 
    satisfies $q(\hat{t}) = q_{\mathrm{ad}}$ for some $\hat{t}>0$, then the right hand derivative of $q$ satisfies for 
    \begin{equation*}
        \frac{d^+q}{dt} (t) = \left\{ \begin{matrix}
            0 &\,, & \alpha(q_{\mathrm{ad}}) R(t) \geq \mu^*,\\
            \alpha(q(t))R(t) - \mu(q(t)) &\,,& \text{else.}
        \end{matrix}\right.
    \end{equation*}
    Let $R_{\mathrm{ad}}> 0$ be the unique point for which $\alpha(q_{\mathrm{ad}})R_{\mathrm{ad}} = \mu^*$, so that $(R_{\mathrm{ad}},q_{\mathrm{ad}})$ is the intersection of the chattering surface with the red line depicting the condition $\dot{q} = 0$ in Figure~\ref{fig:domain2}.
    Note that if $q = q_{\mathrm{ad}}$, $R\geq R_{\mathrm{ad}}$ we have
    \begin{multline*}
        \dot R = K_R - f(p(q_{\mathrm{ad}})) R - \alpha_{\mathrm{ad}}(R,q_{\mathrm{ad}}) \\ = K_R - f(p(q_{\mathrm{ad}})) R - \alpha(q_{\mathrm{ad}}) R_{\mathrm{ad}} \\ \leq K_R - f(p(q_{\mathrm{ad}})) R_{\mathrm{ad}} - \alpha(q_{\mathrm{ad}}) R_{\mathrm{ad}} < 0.
    \end{multline*}
    Thus the trajectory enters the rectangle with vertex points $(0,0)$ and $(R_{\mathrm{ad}},q_{\mathrm{ad}})$ in finite time. By Theorem~\ref{theo:domain} this rectangle is forward invariant under the dynamics of \eqref{eq:system} and contained in the domain of attraction of $x^s$. This shows the assertion.
\end{proof}

\section{Multiple user classes}
\label{sec:multipleusers}

In this section we extend the analysis from a single user queue to the case that there are different types of users which are distinguished by their respective price sensitivity, and where one of these additional classes may correspond to an unresponsive class. As discussed in the introduction this may be a class of users that does not react to price signals, but also different levels of sensitivity can be included in the analysis. From now on we abbreviate $\underline{n}:= \{1,\ldots,n\}$ for $n \in \N$. For the access system there are now $n \in \N$ user classes that all enter the main queue by the same ratio $\alpha(q)$ and which all receive the same price signal $p(q)$. The subclass sensitivity to the price signal of each user group is denoted $f_i$, $i\in \underline{n}$ with $f_n(\cdot)=0$ (the non-responsive class). 
This results in the system
\begin{subequations}
\label{eq:systemRUq}
\begin{align}
\label{eq:systemR11}
    \dot R_i &= K_i - f_i(p(q)) R_i - \alpha(q) R_i, \quad i\in \underline{n},\\
\label{eq:systemq1}
    \dot q &= \alpha(q) \sum_{i=1}^n R_i - \mu(q).
\end{align}
\end{subequations}
\begin{assumption}
\label{ass:pricesensitivities}
    The class price response functions $f_i: \R_+ \to \R_+$, 
    $i \in \underline{n}$,
    are each nondecreasing and continuously differentiable. 
    Furthermore, one of the price response functions has strictly positive derivative.\end{assumption}

Note that Assumption~\ref{ass:pricesensitivities} also admits the important case that one particular user class is unres\-pon\-sive to the price signal, e.g. by setting $f_n(\cdot)=0$.

We begin with the consideration of constant inflow rates $K_i > 0$, $i \in \underline{n}$. Given fixed inflow rates $K_i$, $i=1,\ldots,n$, fixed points $(R^*,q^*)\in \R^{n+1}_+$ of \eqref{eq:systemRUq} are characterized by the conditions
\begin{subequations}
    \label{eq:fpconds-musers}
    \begin{align}
        \label{eq:fpconds-musers-R}
    R_i^* &= \frac{K_i}{f_i(p(q^*)) + \alpha(q^*)}, \quad i= 1,\ldots,n, \\       
        \frac{\mu(q^*)}{\alpha(q^*)} &=  \sum_{i=1}^n R_i^* .
        \label{eq:fpconds-musers-q}
    \end{align}
\end{subequations}
We note that the first set of equations shows that a potential value of $q^*$ uniquely determines the corresponding $R_i^*$, $i\in \underline{n}$ such that $(q^*,R^*)$ is a fixed point.
To determine potential values $q^*$ consider the function
\begin{equation}
\label{eq:g-def}
   g(q):= \sum_{i=1}^n  \frac{ K_i}{f_i(p(q)) + \alpha(q)} - \frac{\mu(q)}{\alpha(q)}, \quad q \geq 0.
\end{equation}
The points $q^*$ for which there exists a solution of \eqref{eq:fpconds-musers} are precisely the zeros of $g$. Note that by assumption $g(0) = (1/\alpha(0)) \sum_{i=1}^n K_i >0$. 
In addition, we will assume the following.

\begin{assumption}
    \label{ass:g-assumption}
 The function $g$ in \eqref{eq:g-def} has exactly two zeros  $0<q^s < q_m < q^u < 2q_m$.  In addition, $g$ undergoes a sign change in each of its zeros and in $q^u$ this sign change is transversal, i.e.,
\begin{equation}
    \label{eq:transverse-cond-musers}
    g'(q^u) = - \sum_{i=1}^n
    \frac{K_i ((f_i\circ p)'(q^u) + \alpha'(q^u))}{(f_i(p(q^u))  + \alpha(q^u))^2 } - \frac{\mu'(q^u)\alpha(q^u) - \mu(q^u) \alpha'(q^u)}{\alpha^2(q^u)} >0.  
\end{equation}   
\end{assumption}

Note that if Assumption~\ref{ass:g-assumption} holds then the values $q^s,q^u$ are characterized by the conditions $p'(q^s) >0$, $p'(q^u) <0$.
While the notation is arbitrary at this point, it will turn out later that indeed $q^s$ gives rise to a unique asymptotically stable fixed point $(R^s,q^s)$, while $q^u$ corresponds to an unstable fixed point $(R^u, q^u)$.

We now turn to the local stability analysis of the two fixed points. To this end, 
the Jacobian of the vector field in its points of differentiability $(R,q) \in \R^{n+1}_+$ is given by (again omitting the argument $q$ for legibility) 
\begin{equation}
\label{eq:DFRq-ndim}
DF(R,q) =     \begin{bmatrix}
        - \vartheta_1 & 0  &\ldots & 0 & 
        - \vartheta_1' R_1 \\
        0 & - \vartheta_2 & \ddots &\vdots  & 
        - \vartheta_2' R_2 \\
        \vdots&\ddots & \ddots & 0 & \vdots \\
         0 & \ldots &0& - \vartheta_n  & - \vartheta_n' R_n \\
        \alpha & & \ldots & \alpha & \alpha' (\sum_{i=1}^n R_i) - \mu' 
    \end{bmatrix}, 
\end{equation}
where $\vartheta_i := (f_i\circ p + \alpha)$, $i\in\underline{n}$. 

Before proceeding to analyze the local stability of the fixed points, we recall the following definitions and results from \cite{kaszkurewicz2012matrix}.
\begin{itemize}
\item[(i)] A matrix $A=(a_{ij})_{i,j=1}^n\in \R^{n\times n}$ is said to be indecomposable (irreducible) if its underlying directed graph is strongly connected, meaning there is a path from every node to every other node in the graph. 

\item[(ii)] The matrix $A$ is said to be {\em combinatorially symmetric} if  $a_{ij} \neq 0$ implies that $a_{ji} \neq 0$ for all $j \neq i$.  

\item[(iii)] A {\em chain of length r} of $A$ is a nonzero product of entries $a_{i_1i_2} \cdot a_{i_2i_3} \cdot .... \cdot a_{i_{r-1}i_r} \cdot a_{i_{r}i_{r+1}}$, where $i_1,i_2,..,i_{r}\in \underline{n}$ are distinct. A chain is said to be a {\em cycle of length $r$} if $i_{r+1} = i_1$. If $A$ does not have cycles of length $r \geq k$ then $A$ is said to be {\em acyclic-$k$}.  
\end{itemize}

Recall furthermore that the principal minors of a matrix $A\in \R^{n\times n}$ are the determinants of the $n$ ``upper-left" submatrices of $A$ obtained by deleting the last $k$ rows and columns or $A$, for $k=0,\ldots, n-1$. We use the convention that $\det A_k$ is the minor of the submatrix of size $k\times k$.

Let $A$ be acyclic-3. Associated to $A$ the matrix $A^+$ is defined by setting $a^+_{ii} := a_{ii}$, i.e. by keeping the diagonal, and by setting
\begin{equation*}
    a_{ij}^+ := \left\{ \begin{matrix}
        a_{ij} & \quad & \text{ if } a_{ij}a_{ji} \geq 0\\
        0 && \text{ if } a_{ij}a_{ji} < 0
    \end{matrix}\right., \quad 1 \leq i,j \leq n, i\neq j.
\end{equation*}

\begin{theorem}(\cite[Theorem~2.2.7]{kaszkurewicz2012matrix})
    \label{thm:kask-bhaya}
    An indecomposable acyclic-3 matrix $A\in \R^{n\times n}$ is diagonally stable if and only if the principal minors $\det (A^+_k)$, $k=1,\ldots,n$, of the associated matrix $A^+$ satisfy
    \begin{equation}
    \label{eq:thm-kasj-bhaya}
      (-1)^k \det (A^+_k) >0, \quad k\in\underline{n}.
    \end{equation}
\end{theorem}

Then, the following holds directly from the special structure of $DF(R,q)$ and Theorem \ref{thm:kask-bhaya}.
\begin{theorem}
\label{thm:musers-localstability}
Let Assumptions~\ref{ass:pricesensitivities} and \ref{ass:g-assumption} hold.
    For system \eqref{eq:systemRUq} the fixed point $(R^s,q^s)$ characterized by $p'(q^s) >0$ is locally asymptotically stable.
\end{theorem}

\begin{proof}
We will show that Theorem~\ref{thm:kask-bhaya} is applicable to $A:= DF(R^s,q^s)\in \R^{(n+1)\times (n+1)}$. Then a diagonal quadratic Lyapunov function exists for $A$ and the claim follows from Lyapunov's linearization theorem, \cite{khalil2002nonlinear}. 

We may assume without loss of generality that for all $i\in \underline{n}$ the generic condition $(f_i\circ p)'(q^s) + \alpha'(q^s) \neq 0$ holds.
If this is not the case, we may permute rows and columns so that all zero entries of the last column are at the beginning of the column. Then the lower triangular block structure of the resulting matrix shows that it is sufficient to prove that the maximal lower right block is asymptotically stable in which all elements of the last column are nonzero. 
In this case all entries of the the last row and column of $A$ are nonzero. This implies that $A$ is indecomposable, has exactly $n$ non-zero elements above the diagonal, and is combinatorially symmetric. By \cite[Lemma~2.2.3]{kaszkurewicz2012matrix} this implies that $A$ is acyclic-3.

As the upper-left $n\times n$ submatrix of $A$ is diagonal with negative entries, it only remains to check the determinant of $A^+$ in order to satisfy the conditions of Theorem~\ref{thm:kask-bhaya}. Again we may assume without loss of generality that $a_{i,n+1}^+>0$, $i\in \underline{n}$. Otherwise a permutation reduces the problem to a lower block matrix for which this is the case.

Using the Schur complement formula (A.22, page 656 \cite{Kailath}) we obtain
\begin{equation}
\label{eq:det-comp}
   \det A^+ = \left((-1)^n \prod_{k=1}^n (f_k(p) + \alpha) \right) \left( \alpha' \left(\sum_{k=1}^n R_k\right) - \mu'   -\alpha \sum_{k=1}^n \frac{ ( (f_k \circ p)' + \alpha')}{f_k(p) + \alpha}R_k\right).
\end{equation}
To determine the sign of $\det A^+$ it thus just remains to examine the sign of the second factor and in order that \eqref{eq:thm-kasj-bhaya} holds this factor needs to be negative. Rearranging the terms, the expression for the last factor reads
\begin{equation*}
    -\mu' - \alpha \sum_{k=1}^n \frac{ (f_k \circ p)' }{f_k(p) + \alpha}R_k + \alpha' \sum_{k=1}^n \frac{  f_k(p)}{f_k(p) + \alpha}R_k.
\end{equation*}

With our assumptions on the systems parameters, namely $\mu' \geq 0$, $\alpha(q^s)>0, f_k(p(q^s))\geq 0$, $(f_k \circ p)'(q^s) \geq 0$, $\alpha'(q^s) \leq 0$, we see that every summand is nonpositive and the second is strictly negative as for some $k$ we have $(f_k \circ p)'(q^s) > 0$ by Assumption~\ref{ass:pricesensitivities}.

In summary, we have
\begin{equation*}
    (-1)^{n+1} \det A^+ = (-1)^{2(n+1)} |\det A^+| >0.
\end{equation*}
This shows that Theorem~\ref{thm:kask-bhaya} is applicable and there exists a positive definite diagonal Lyapunov function for $A= DF(R^s,q^s)$. 
\end{proof}

The next result shows that the fixed point $(R^u,q^u)$ is unstable.

\begin{theorem}
\label{thm:musers-localinstab}
  Let Assumption~\ref{ass:pricesensitivities} and \ref{ass:g-assumption} hold.
    For system \eqref{eq:systemRUq} consider the fixed point $(R^u,q^u)$ characterized by $p'(q^u) <0$.
    The Jacobian $DF(R^u,q^u)$ has the following properties:
    \begin{enumerate}
        \item[(i)] $DF(R^u,q^u)$ has at least $n-1$ real, negative eigenvalues.
        \item[(ii)] $DF(R^u,q^u)$ has a positive eigenvalue.
    \end{enumerate}
    In particular, $(R^u,q^u)$ is an exponentially  unstable fixed point. 
\end{theorem}

\begin{proof}
    (i) As $p'(q^u)<0$, all of the terms $a_{i,n+1}:=-((f_i\circ p)'(q^u) + \alpha'(q^u))R_i^u$, $i\in \underline{n}$, appearing in the last column of $A:=DF(R^u,q^u)$ are positive. We may thus apply a diagonal similarity transformation $T= \operatorname{diag}(t_1,\ldots,t_n,1)$, with $t_i = \sqrt{\alpha/a_{i,n+1}}$, $i\in \underline{n}$, and obtain that $TAT^{-1}$ is symmetric. As the diagonal part $\operatorname{diag}(a_{11},\ldots,a_{nn})$ is a principal submatrix of $TAT^{-1}$ it follows from Cauchy's interlacing theorem, \cite[Theorem~IV.4.2]{StewSun90}, that the eigenvalues $\lambda_1\geq \lambda_2 \geq \dots \geq \lambda_{n+1}$ of $A$ (which coincide with the eigenvalues of $TAT^{-1}$) satisfy the interlacing property
    \begin{equation*}
        \lambda_1 \geq \mu_1 \geq \lambda_2 \geq \mu_2 \geq \dots \geq \mu_n \geq \lambda_n, 
    \end{equation*}
    where the $\mu_1,\ldots,\mu_n$ are the ordered list of the entries $a_{11},\ldots,a_{nn}$. As all the $\mu_i$ are negative, at most the eigenvalue $\lambda_1$ can be positive. 
    
    (ii) By (i), $\lambda_1 >0$ if and only if $(-1)^{n}\det(A) >0$. As in \eqref{eq:det-comp}, this condition holds if and only if the factor
    \begin{equation}
          \alpha' \left(\sum_{k=1}^n R_k\right) - \mu'   -\alpha \sum_{k=1}^n \frac{ ( (f_k \circ p)' + \alpha')}{f_k(p) + \alpha}R_k >0.
    \end{equation}
    Dividing by $\alpha$, substituting $R_k$ by the fixed point condition \eqref{eq:fpconds-musers-R} on the right and $\sum_{k=1}^n R_k$ by the fixed point condition \eqref{eq:fpconds-musers-q} on the left, we see that this condition is equivalent to \eqref{eq:transverse-cond-musers}. It thus follows from our assumptions that $\lambda_1>0$. Then by linearization, \cite{khalil2002nonlinear}, the fixed point $(R^u,q^u)$ of system \eqref{eq:systemRUq} is exponentially unstable. This completes the proof.
\end{proof}

\subsection{Invariance}
\label{subsec:imvariance-musers}

In this section we study forward invariance properties of the system \eqref{eq:systemRUq}. To this end, we make the following basic observations for arbitrary points $(R,q) \in \R^n_+ \times \R_+$:
\begin{align}
    \label{eq:Ridot-neq}
    \dot R_i < 0 & \quad \Leftrightarrow \quad R_i > \frac{K_i}{\alpha(q) + f_i(p(q))}, \quad i \in \underline{n},\\
    \dot q < 0 & \quad \Leftrightarrow \quad \sum_{i=1}^n R_i < \frac{\mu(q)}{\alpha(q)}.
    \label{eq:qdot-neg}
\end{align}

In order to study invariance properties we make a further assumption.

\begin{assumption}
    \label{ass:forinvariance}
 Assume that for all $i\in \underline{n}$ we have
    \begin{equation}
    \label{eq:minv-minimality}
      f_i(p(q^u)) + \alpha(q^u) = \min_{q\in[0,q^u]} f_i(p(q)) + \alpha(q).
    \end{equation}   
\end{assumption}

On the interval $[q_m,q_u]$ the functions $f_i\circ p(q) + \alpha(q)$, $i\in \underline{n}$, are all strictly decreasing. Similarly to the discussion in \eqref{eq:Rmaxdef} we define $q^\dagger$ as the minimal value in $[q_m,q_u]$ such that for all $\hat{q} \in [q^\dagger, q^u]$ we have for all $i\in \underline{n}$
    \begin{equation}
    \label{eq:min-minimality-hat}
      f_i(p(\hat{q})) + \alpha(\hat{q}) = \min_{q\in[0,\hat{q}]} f_i(p(q)) + \alpha(q).
    \end{equation}   
Given $\hat{q} \in [q^\dagger, q^u]$, define
\begin{equation}
\label{eq:Rhat-def}
    R_i(\hat{q}) = \frac{K_i}{f_i(p(\hat{q}))+ \alpha(\hat{q})}
\end{equation}
and the associate "box" in $\R_+^{n+1}$ denoted $B(\hat{q})$ and defined by
\begin{equation}
    \label{eq:Boxdef}
     \{ (R,q) \in \R^{n+1}_+ \; \vert \; 0 \leq q \leq \hat{q}, \\ 0\leq R_i \leq R_i(\hat{q}), i\in \underline{n}\}.
\end{equation}

\begin{proposition}
\label{lem:n-users-invariant-box}
    Let Assumptions~\ref{ass:pricesensitivities}, \ref{ass:g-assumption}, and \ref{ass:forinvariance} hold.
    Let $q^s, q^{u}$ be the two solutions of \eqref{eq:g-def} and $q^\dagger$ as defined before \eqref{eq:min-minimality-hat}.
For every $\hat{q}\in [q^\dagger, q^u]$ the box $B(\hat{q})$
is forward invariant under the flow of the system \eqref{eq:systemRUq}.
\end{proposition}

\begin{proof}
    The box $B(\hat{q})$ is a polytope and it is again sufficient to consider
    the outward normals on the facets characterized by $R_i = 0$, or $R_i= \hat{R}_i:=R_i(\hat{q})$, respectively $q=0$, $q=\hat{q}$. These are the respective standard unit vectors $-e_i$, or $e_i$. The invariance conditions are thus for $i\in \underline{n}$
    \begin{subequations}
        \begin{align}
            \label{eq:minvariance-Ri0}
            \langle -e_i, F(R,q) \rangle &= - K_i \leq 0, &&\text{if } (R,q) \in B(\hat{q}), R_i = 0; 
        \\
        \label{eq:minvariance-Riu}
            \langle e_i, F(R,q) \rangle &= K_i - f_i(p(q))\hat{R}_i - \alpha(q)\hat{R}_i \leq 0,  &&\text{if }  (R,q) \in B(\hat{q}), R_i = \hat{R}_i; 
\intertext{and for $i=n+1$}
           \label{eq:minvariance-q0}
        \langle -e_{n+1}, F(R,q) \rangle &= 
        - \alpha(0) \sum_{i=1}^n R_i \leq 0,
          &&\text{if }  (R,q) \in B(\hat{q}),q = 0;
\\
           \label{eq:minvariance-qu}
       \langle e_{n+1}, F(R,q) \rangle &= \alpha(\hat{q}) \sum_{i=1}^n R_i - \mu(\hat{q})\leq 0,
           &&\text{if }  (R,q) \in B(\hat{q}), q = \hat{q}.
\end{align}
    \end{subequations}
Conditions \eqref{eq:minvariance-Ri0}, \eqref{eq:minvariance-q0} are trivially satisfied. Condition \eqref{eq:minvariance-Riu} is a consequence of \eqref{eq:min-minimality-hat} and \eqref{eq:minvariance-qu} follows from a combination of \eqref{eq:Rhat-def} and \eqref{eq:g-def} together with Assumption~\ref{ass:g-assumption} by which $g(\qad) <0$.
The forward invariance of $B(\hat{q})$ now follows from the invariance principle, \cite[Theorem 4.3.8]{ClarLedy98}.
\end{proof}

\begin{corollary}
\label{cor:invattractor}
Under the conditions of Proposition~\ref{lem:n-users-invariant-box}, for every $\hat{q}\in [q^\dagger,q^u]$ the flow on
$B(\hat{q})$ is volume contracting. In particular, there exists a compact
invariant set $C(\hat{q}) \subset B(\hat{q})$. The invariant set $C(\hat{q})$ attracts all trajectories with initial condition  
$x \in B(\hat{q})$ in the sense that
\[
\lim_{t\to \infty} \dist(\varphi(t,x),C(\hat{q})) = 0.
\]
Moreover, $C(\hat{q})$ has $n$-dimensional Lebesgue measure $0$ and hence has empty interior. For $q^\dagger\leq \hat{q} < q^u$ the set $C(\hat{q})$ is independent of $\hat{q}$ and equals $C(q^\dagger)$.
\end{corollary}

\begin{proof}
    On $B(\hat{q})$ the divergence, i.e. the trace of $DF(R,q)$, is strictly negative, see \eqref{eq:DFRq-ndim}. It follows as before that the flow on $B(\hat{q})$ is strictly volume contracting and hence
    \begin{equation}
        \lim_{t\to \infty} \mathrm{vol}_n(\varphi(t,B(\hat{q}))) =0.
    \end{equation}
    In addition, by the forward invariance established in Proposition~\ref{lem:n-users-invariant-box} and continuity of the flow, the family $(\varphi(t,B(\hat{q})))_{t\geq 0}$ is a decreasing family of compact sets. Thus
    \begin{equation}
        C(\hat{q}) := \bigcap_{t\geq 0} \varphi(t,B(\hat{q}))
    \end{equation}
    is by construction a nonempty, compact, and invariant set that attracts all trajectories starting in $B(\hat{q})$. By the previous consideration the $n$-dimensional Lebesgue measure of $C(\hat{q})$ is $0$.

    Finally, if $q^\dagger < \hat{q} < q^u$, then \eqref{eq:min-minimality-hat} shows that $f_i(p(\hat{q}))+ \alpha(\hat{q})$ is a strict minimum over the interval $[0,\hat{q}]$. This means that equality to zero in \eqref{eq:minvariance-Riu} is only attained for $q=\hat{q}$. Similarly, equality to zero in \eqref{eq:minvariance-qu} is not attained as $g(\hat{q})< 0$. As $(\hat{R},\hat{q})$ is not a fixed point, this means that $C(\hat{q})$ cannot intersect the facets of $B(\hat{q})$ defined by $q=\hat{q}$ or $R_i = R_i(\hat{q}))$. By compactness of $C(\hat{q})$ there thus exists a $q^\dagger \leq q_2 < \hat{q}$ with $C(\hat{q}) \subset B(q_2)$ and in particular $C(\hat{q})=C(q_2)$. As $\hat{q}$ was arbitrary in the open interval $(q^\dagger, q^u)$, this shows that $C(\hat{q}) = C(q^\dagger)$ for all $\hat{q}\in (q^\dagger, q^u)$.
\end{proof}
Note that the final argument of the previous proof does not apply to $B(q^u)$ as $(R^u,q^u)$ is a fixed point.\newline

\subsection{Modified admission policy}
As in Subsection~\ref{subsec:modified-1user} we now introduce a chattering surface at a value $q^\dagger < q_{ad} <q^u$ at which the admission to the service queue is reduced to the current service capacity. The difference to the expression for $\alphad$ in \eqref{eq:chattering_queue_admittance} is merely that the admittance to the queue is now a function of the sum $\sum_{i=1}^n R_i$. There are various options how the available admittance to the queue $\alpha_{ad}(\sum_{i=1}^n R_i,q)$ is distributed among the different user classes. Here we take the approach that this distribution is proportionate to queue length. This models the common practical case that user classes cannot be distinguished by the admittance process to the service queue. With this modeling assumption the extension of \eqref{eq:systemRUq} is 
\begin{subequations}
\label{eq:system-chatter-nusuers}
\begin{align}
\label{eq:systemR-chatt-nusers}
    \dot R_i &= K_i - f_i(p(q)) R_i - \alpha_{\mathrm{ad}}\left(\sum_{k=1}^n R_k,q\right)\frac{R_i}{\sum_{k=1}^n R_k},\\
\label{eq:systemq-chatt-nusers}
    \dot q &= \alpha_{\mathrm{ad}}\left(\sum_{k=1}^n R_k,q\right) - \mu(q).
\end{align}
\end{subequations}

If $q<q_{\mathrm{ad}}$, the equations in \eqref{eq:system-chatter-nusuers} coincide with those of \eqref{eq:systemRUq}, as in this case $\alpha_{\mathrm{ad}}\left(\sum_{k=1}^n R_k,q\right) = \alpha(q)\sum_{k=1}^n R_k$.
For the following result we continue to use the notation introduced in \eqref{eq:min-minimality-hat}-\eqref{eq:Boxdef}.

\begin{boxedminipage}{0.98\textwidth} 
\begin{theorem}{\bf Main Result (multiple responsive classes)}
\label{th:main_result_sec4}
Let Assumptions~\ref{ass:pricesensitivities}, \ref{ass:g-assumption}, and \ref{ass:forinvariance} hold and let $\qad \in (q^\dagger, q^u)$.
    Consider \eqref{eq:system-chatter-nusuers} with \eqref{eq:chattering_queue_admittance} and a continuously differentiable admission rate $\alpha(q)$.
    Then $x^s=(R^s,q^s)$ is a locally asymptotically stable fixed point of \eqref{eq:system-chatter-nusuers} and for every initial condition $x=(R,q) \in X_{\mathrm{ad}}:=[0,\infty)^{n} \times [0,q_{\mathrm{ad}}]$ we have for the trajectory $\varphi(\cdot;x)$
    \begin{equation}
    \label{eq:Bqad-limit}
        \lim_{t\to \infty} \dist(\varphi(t;x),C(q^\dagger)) = 0,
    \end{equation}
    with $C(q^\dagger)$ given by Corollary~\ref{cor:invattractor}.
\end{theorem}
\end{boxedminipage}

\begin{proof}
For $\qad \in  [q^\dagger, q^u)$ consider the box $B(\qad)$ defined in \eqref{eq:Boxdef}. By Assumption~\ref{ass:g-assumption} we have $g(\qad)<0$ and this implies for $\hat{R}_i:=R_i(\qad)$, $i\in\underline{n}$, defined in \eqref{eq:Rhat-def} that 
\begin{equation}
\label{eq:map-1}
    \alpha(\qad) \sum_{i=1}^n R_i(\qad) < \mu^*.
\end{equation}
This implies that on $B(\qad)$ the dynamics of \eqref{eq:systemRUq} and \eqref{eq:system-chatter-nusuers} coincide. By Proposition~\ref{prop:invariance} the box $B(\qad)$ is forward invariant. 

We will show that on $X_{\mathrm{ad}} \setminus B(\qad)$ the function $V: X_{\mathrm{ad}} \to \R_+$ defined by
\begin{equation}
    V(R,q) = \max_{i\in \underline{n}}\{ R_i/\hat{R}_i\}
\end{equation} 
is strictly decreasing along trajectories of \eqref{eq:system-chatter-nusuers}. To this end, let $(R,q)\in X_{\mathrm{ad}} \setminus B(\qad)$ and let $i\in \underline{n}$ satisfy
\begin{equation}
\label{eq:maxRi}
    \frac{R_i}{\hat{R}_i} = \max_{j\in \underline{n}} \frac{R_j}{\hat{R}_j}.
\end{equation}
We distinguish two cases: (i) if $q<\qad$ or $q=\qad$ and
$\alpha(\qad) \sum_{j=1}^n R_j < \mu^*$, then the dynamics of \eqref{eq:system-chatter-nusuers} are given by \eqref{eq:systemRUq} and we have
\begin{multline*}
    \dot R_i = K_i - f_i(p(q)) R_i - \alpha(q) R_i 
    \stackrel{\eqref{eq:min-minimality-hat}}{<} K_i - f_i(p(\qad)) R_i - \alpha(\qad) R_i \\
    \stackrel{R_i > \hat{R}_i}{<} K_i - f_i(p(\qad)) \hat{R}_i - \alpha(\qad) \hat{R}_i \stackrel{\eqref{eq:Rhat-def}}{=} 0.
\end{multline*}
As this holds for all indices $i$ satisfying \eqref{eq:maxRi} this shows that $V$ is decreasing in $(R,q)$.

(ii) Now consider $q=\qad$ and $\alpha(\qad)\sum_{j=1}^n R_j = \mu^* $. Then we have in the point $(R,\qad)$ for \eqref{eq:system-chatter-nusuers}
\begin{align*}
    \dot R_i &= K_i - f_i(p(\qad)) R_i - \alphad\left(\sum_{j=1}^n R_j, \qad\right) \frac{R_i}{\sum_{j=1}^n R_j}\\
    &= K_i - f_i(p(\qad)) R_i - \mu^* \frac{R_i}{\sum_{j=1}^n R_j}\\
    &= K_i - f_i(p(\qad)) R_i - \alpha(\qad)R_i \\
    &< K_i - f_i(p(\qad)) \hat{R}_i - \alpha(\qad)\hat{R}_i=0.
\end{align*}
For the general case $\sum_{j=1}^n R_j > \mu^* / \alpha(\qad)$, choose $\lambda >0$ such that $\lambda\sum_{j=1}^n R_j = \mu^* / \alpha(\qad)$. This does not affect the maximizing indices $i$ in the definition of $V$. By the previous computation we have 
\begin{align*}
    \dot R_i &= K_i - f_i(p(\qad)) R_i - \mu^* \frac{R_i}{\sum_{j=1}^n R_j}\\
    &\leq K_i - f_i(p(\qad)) \lambda R_i - \mu^* \frac{\lambda R_i}{\sum_{j=1}^n \lambda R_j} <0.
\end{align*}
Again, $V$ is decreasing in $(R,\qad)$.

Note that $V(R,q) > 1$ if and only if $(R,q) \in X_{\mathrm{ad}} \setminus B(\qad)$ and so the strict decrease of $V$ along trajectories shows that for all $x \in X_{\mathrm{ad}}$ we have $\lim_{t\to\infty} \dist
(\varphi(t,x),B(\qad)) =0$. Take $\varepsilon>0$ sufficiently small such that $\qad + \varepsilon < q^u$. Then the previous argument shows that for all $x \in X_{\mathrm{ad}}$ there exists a $T>0$ such that 
$\varphi(T,x) \in B(\qad + \varepsilon)$. By Corollary~\ref{cor:invattractor} applied to $B(\qad + \varepsilon)$ we obtain \eqref{eq:Bqad-limit}.
\end{proof}

\section{Interpretation of main results} 
\label{sec:interpretation}
Theorem \ref{th:main_result_sec3} and Theorem \ref{th:main_result_sec4} constitute the main technical results of the paper. Theorem \ref{th:main_result_sec3} states that for any initial condition, the unforced system will return to the stable equilibrium $x^s$ contained in the interior of the invariant set. Theorem \ref{th:main_result_sec4} states that in the case of multiple user classes, the unforced system will return to the invariant set that contains the stable equilibrium point.\newline 

To appreciate the importance of these results, consider briefly the situation when elastic and inelastic traffic compete (Equations \eqref{eq:system_bob}) which we have written in a slightly different form:
\begin{subequations}
\label{eq:systemR11_bob2}
\begin{align}
\label{eq:systemR12_bob2}
    \dot R &= K_R - f(p(q)) R - \zeta_1(R,U,q)R, \\
\label{eq:systemU_bob2}
    \dot U &= K_U(t) - \zeta_2(R,U,q) U, \\
\label{eq:systemq1_bob2}
\dot q &= \alphad(R+U, q) - \mu(q).
\end{align}
\end{subequations}
where for all $q < \qad$, $\zeta_1(R,U,q) = \zeta_2(R,U,q) = \alpha(q)$ and for all $q \geq \qad$, $\zeta_1(R,U,q)U+\zeta_2(R,U,q)R = \alphad(R+U,q)$. In addition, we assume that $\zeta_1,\zeta_2$ are continuous, strictly positive on the set $\{(R,U,q); R,U\geq 0, q \geq q_{\mathrm{ad}}\}$,
and for all bounds $K>0$ there exists a constant $c>0$ such that $q,R\leq K$ implies $\zeta_2(R,U,q)>c$ for all $U>0$.
The functions $\zeta_i$ constrain admission to the service queue so that conservation of flow is satisfied. When the 
non-responsive flow is present (which acts like an input) we expect that the queue will tend to fill, both due to reduced price for the elastic traffic and 
and due to the non-responsive flow which does not respond to price at all. When the non-responsive traffic leaves the system we recover   
the dynamics \eqref{eq:system-chatter}. 
The main consequence of Theorem \ref{th:main_result_sec3} is that these dynamics will stabilize at the stable equilibrium point and 
this will happen without the need to observe the departure of the non-elastic traffic.
The following result characterizes the behavior of the admission system in case of a surge of inelastic demand.
\begin{theorem}
    Consider system~\eqref{eq:systemR11_bob2} with $\alpha_{\mathrm{ad}}$ given by \eqref{eq:chattering_queue_admittance}. Let $K_U:\R_+\to \R_+$ be continuous and satisfy
    \begin{equation}
        \lim_{t\to \infty} K_U(t) = 0.
    \end{equation}
    Then for all initial conditions $x_0=(R_0,U_0,q_0)\in R_+^2\times [0,q_{ad}]$ 
    \begin{equation}
        \label{eq:surge-limit}
        \lim \varphi(t,x_0) = (R^s,0,q^s),
    \end{equation}
    where $(R^s,q^s)$ is the stable fixed point of system \eqref{eq:system-chatter}.
\end{theorem}

\begin{proof}
Fix the initial condition $x_0$. We will first show that the corresponding trajectory is uniformly bounded. This is clear from \eqref{eq:systemU_bob2} for the $q$-component. The $R$-component remains bounded, as $\dot R <0$ if $q=\qad$ and $R> K_R / f(p(\qad))$ and on the other hand if $R(t) > \alpha(\qad) / \mu^*$ on an interval of the form $[T,\infty)$, then $q(t) = \qad$ for all $t$ sufficiently large.
Finally, by assumption there is a constant $c >0$ such $\zeta_2(q(t),U(t),R(t)) > c$ for all $t\geq 0$. Then \eqref{eq:systemU_bob2} implies using the variation of constants formula that 
\begin{equation}
    0\leq U(t) \leq e^{-ct}U_0 + e^{-ct}\int_0^t e^{cs} K_U(s) ds.
\end{equation}
The expression on the right tends to $0$ for $t\to\infty$ by \cite[Lemma 8.3.18]{hinrichsen2026mathematicalII}. As $U(t) \to 0$ for $t\to \infty$, the $(R,q)$-dynamics of \eqref{eq:systemR11_bob2} are asymptotically given by \eqref{eq:system-chatter}. As the $(R,q)$-component of the trajectory is bounded, it thus converges to $(R^s,q^s)$ by Theorem~\ref{th:main_result_sec3}.
\end{proof}

\begin{remark}
   A similar result may be obtained in the more general multi-class case by invoking Theorem \ref{th:main_result_sec4}. In this case the conclusion is weaker in that the system will return to the invariant set $C(q^\dagger)$ that contains the stable equilibrium, rather than the equilibrium itself, once the unresponsive traffic has left the system. 
\end{remark}

\section{Non-vanishing prices}
\label{nozeroprice}

We now conclude the analytical part of the paper by modifying the model by introducing a positive saturation in the price function $p(\cdot)$.  As we shall see this yields an alternative to the chattering solution that we have discussed. In addition, we also account for the influx of unresponsive load to the model. For the moment, this is taken to be a deterministic influx $U$, possibly time-dependent. First we modify the definition of the price $p(q)$ from \eqref{eq:costfunction} and set for some $q_n \in (q_m, 2q_m)$
\begin{equation}
\label{eq:costfunction-sat}
    p_{\mathrm{sat}}(q) = \left\{ \,\begin{matrix}
    \beta q & ,&  0\leq q \leq q_m \\
    \beta ( 2 q_m -q) & ,& q_m\leq q \leq q_n  \\
    \beta(2 q_m - q_n) &,& q_n \leq q < \infty \\    
    \end{matrix} \right..
\end{equation}
In addition, we will assume for the moment, that there is a constant load $U$ into the queue represented by $q$. The equations for system \eqref{eq:system} are then modified to
\begin{subequations}
\label{eq:system-sat}
\begin{align}
\label{eq:system-satR}
    \dot R &= K_R - f(p_{\mathrm{sat}}(q)) R - \alpha(q) R,\\
\label{eq:system-satq}
    \dot q &= \alpha(q) R - \mu(q) + K_U.
\end{align}
\end{subequations}

The conditions for fixed points are now, similarly to \eqref{eq:fp-conds3}, of the form
\begin{align}
\begin{split}
    &R^* = \frac{\mu(q^*)-K_U}{\alpha(q^*)}, \\
    &\frac{K_R-\mu(q^*) + K_U}{\mu(q^*) - K_U} = \frac{f(p_{\mathrm{sat}}(q^*))}{\alpha(q^*)}.
\end{split}
\end{align}

In particular, it is necessary, that $\mu^* > K_U$ so that fixed points can exist in $\R^2_+$.
In contrast to the previous section we have enforced new forward invariance properties. As the divergence is still negative everywhere on $\R_+^2$ we also obtain a global result for $\omega$-limit sets. Note that, in contrast system \eqref{eq:system}, there now exists a constant $c_{\mathrm{sat}}>0$ such that
\begin{equation*}
    \alpha(q) + f(p_{\mathrm{sat}}(q)) > c_{\mathrm{sat}} >0 , \quad q \geq 0.
\end{equation*}

\begin{proposition}
    \label{prop:sat2d-omega-limit}
    Consider system \eqref{eq:system-sat} with price function given by \eqref{eq:costfunction-sat}. Assume that $\mu^* > K_U$. Then
    \begin{enumerate}[(i)]
        \item The sets $\R_2^+$, $[0,\infty)\times [0,q_{\max}]$, $[0, \frac{K_U}{c_{\mathrm{sat}}}]\times [0,q_{\max}]$ are forward invariant.
        \item For every $x_0 \in \R_+^2$ we have $\omega(x_0) \subset [0, \frac{K_U}{c_{\mathrm{sat}}} ]\times [0 ,q_{\max}]$ and all $\omega$-limit sets are fixed points. 
    \end{enumerate}
\end{proposition}

\begin{proof} The claims follow by a combination of the arguments presented in the proofs of  Lemma~\ref{lemma:localstability-2d} and Theorem~\ref{theo:domain} together with the fact that $\dot R \leq  K_U - c_{\mathrm{sat}}R < 0$, provided that $R> \frac{K_U}{c_{\mathrm{sat}}}$. 
\end{proof}

\begin{figure}[t]
    \centering
    \includegraphics[width=0.7\linewidth]{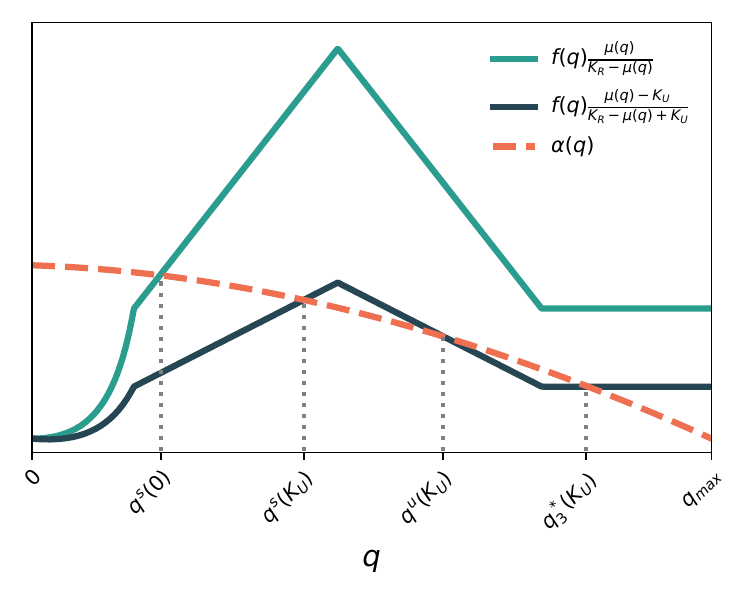}
    \caption{Sketch of the fixed points in the case of the saturated price function of Eq.~\eqref{eq:costfunction-sat}.}
    \label{fig:fixedpoints-saturated}
\end{figure}

A comparison of the new fixed point(s) when unresponsive traffic is and is not present (for the saturated case) is shown in Fig.~\ref{fig:fixedpoints-saturated}. There are now two distinct scenarios depending on the parameters. First, there may be a unique fixed point, which is the globally asymptotically stable with respect to the invariant set $\R_+^2$, or, second, there are three fixed points. In the latter case, two of these points retain the properties of the fixed points $x^s(K_U),x^u(K_U)$ discussed for system \eqref{eq:system}, while a third asymptotically stable fixed point $x_3^*(K_U)$ exists in the large queue size regime. 
The interest of this third fixed point is 
that if the constant influx $K_U$ is temporary and switches back to a lower value $K_U'$ (possibly zero), then the fixed point $x^s(K_U')$ attracts the prior fixed point $x_3^*(K_U)$
In this case, however, one must select the system parameters to avoid a stabilizing solution when the unresponsive traffic leaves the system. It is for this reason that the chattering approach may be viewed as a more robust design.

\section{Simulation Results}
\label{sec:sim}

In this section we use
simulations to implement\footnote{All numerical numerical simulations are implemented in Python with the code available at 
\url{https://github.com/h0mayoun/access-pricing}} and evaluate the efficacy of the proposed pricing control mechanisms. The simulations implement the model and pricing mechanisms numerically using a standard ODE solver applied over time. In particular, we compare the proposed mechanism to a standard surge pricing policy. 
Consistent with the overall system we have presented,
the implemented scenario involves price-responsive ($R$) and unresponsive ($U$) users arriving to receive service with a service queue capacity of $q_{max}=100$ and service rate defined as in Eq.~$\eqref{eq:mudef}$, where we also set $q_c=35<q_{max}$. The queuing users receive service based on an admission function $\alpha(q)$ (dynamically changing with the service queue occupancy $q$) which we take to have the same form as in Fig.~\ref{fig:fixedpoints-saturated} (see the dashed red curve) tuned for the system to have the desired either triple or single 
fixed points in the presence (absence) of the unresponsive group of users when the price is determined by the saturated price function of Eq.~\eqref{eq:costfunction-sat}. Specifically, a third-order polynomial $\alpha:q\rightarrow \R$ is used with the coefficients set for $\alpha(q)$ to be monotonically decreasing over $q\in[0, q_{max}]$ and to mimic the scenario of Fig.~\ref{fig:fixedpoints-saturated} in terms of the fixed points. 

We run the simulation of a scenario that demonstrates different behaviors of the dynamic system: convergence in the presence of only price-responsive users, instability induced by the proposed pricing in the scenario when unresponsive traffic is present, and the bounceback of the system after the termination of the unresponsive demand. The numerical simulations are initialized with $R(0)=50$ and $q(0)=15$ and constant $K_R=4$ (with no unresponsive user arrival $K_U=0$) over $t\in[0, 100]$. The unresponsive demand initiates at $t=100$ and terminates at $t=300$ with $K_U=4$ being constant over this period. At $t=300$, the unresponsive demand terminates and the simulation continues until $t=400$. 

\begin{figure}[htp]
    \centering
    \includegraphics[width=0.65\linewidth]{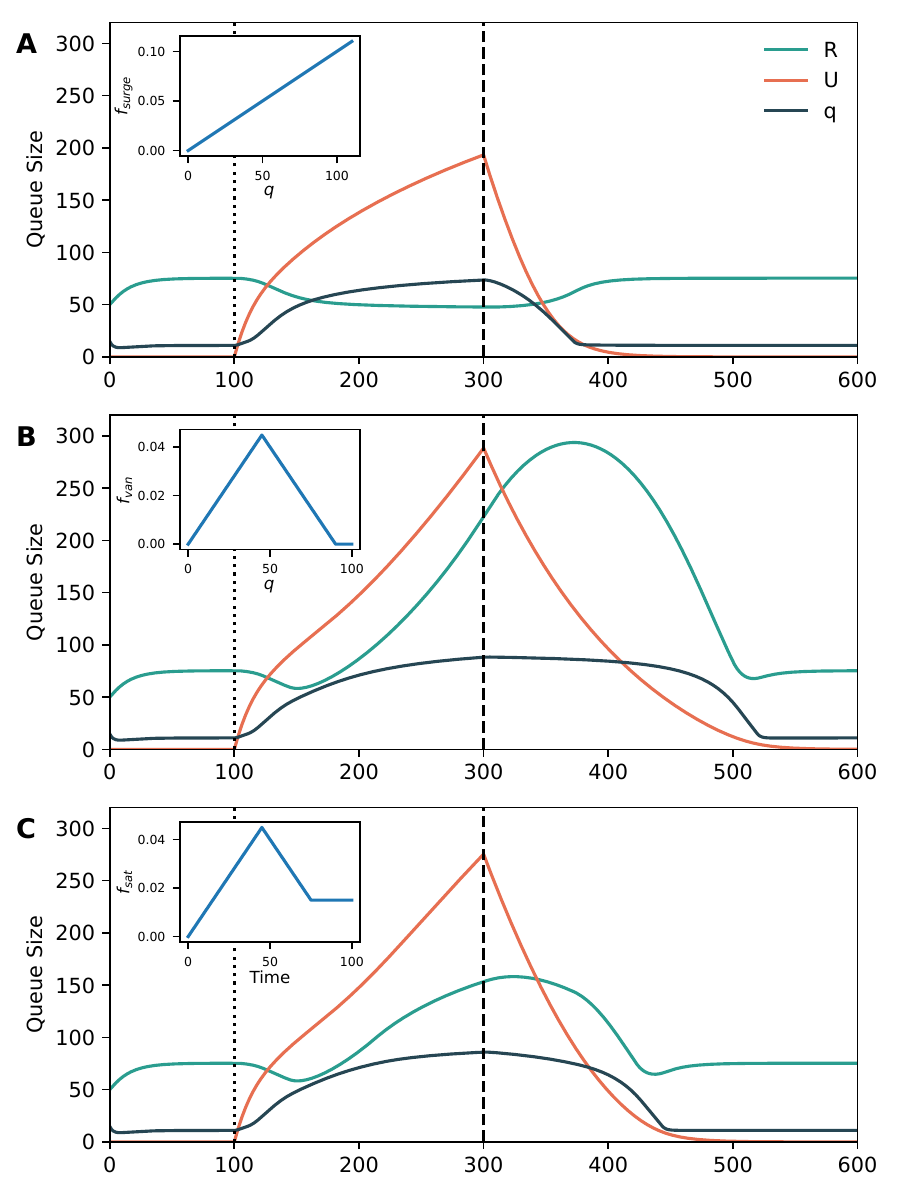}
    \caption{Time series of unresponsive $U$ (red), responsive $R$ (green), and service $q$ queues, compared between three scenarios where  different pricing mechanisms are governing the system: the standard surge pricing (panel A), the proposed mechanism with vanishing price (panel B), and the proposed saturated pricing (panel C). The inset of each panel shows the corresponding pricing function in place. The vertical dotted and dashed lines mark the start and end of the period where unresponsive users enter the system.
    }
    \label{fig:sim_ts}
\end{figure}

We separately simulate the above scenario, first using the standard surge pricing $f_{surge}(q)=\beta q$ (see Fig.~\ref{fig:sim_ts}A), and then demonstrating the proposed pricing mechanism using vanishing price function $f_{van}$ in Eq.~\eqref{eq:costfunction} and saturated price function $f_sat$ in Eq.~\eqref{eq:costfunction-sat} (see Fig.~\ref{fig:sim_ts}B and C, respectively); setting $\beta=10^{-3}$ for both pricing functions, and $q_m=45$, and $q_n=75$ for the proposed mechanism. Figure~\ref{fig:sim_ts} shows the evolving number of responsive ($R$) and unresponsive ($U$) users waiting for the service and the size of the service queue ($q$) over time $t$, comparing the standard versus the proposed pricing mechanism. The inset of each panel in Fig.~\ref{fig:sim_ts} shows the pricing function corresponding to the generated time series of the queue size for the three main compartments of the system. 

Figure~\ref{fig:sim_ts} visualizes the price-responsive population being priced out of the system as a consequence of standard surge pricing over $t\in[100, 300]$, as expected. The proposed pricing only increases the price up to a certain level of service-queue size ($q_m=45$).  When increasing the price in this region does not balance supply and demand, then the price is decreased to balance the competition between price-responsive and unresponsive populations. 
This is seen in how $R$ is reduced in Fig.~\ref{fig:sim_ts}A compared with the size of $R$ increasing at nearly the same rate as $U$ in Figs.~\ref{fig:sim_ts}B,C. By design, the proposed pricing mechanism in Eqs.~\eqref{eq:costfunction-sat} and~\eqref{eq:costfunction} avoids disproportionately pricing out the responsive population to decongest the system, and yet it automatically enables bounce-back to the uncongested phase after a surge of unresponsive user demand. This demonstrates that in contrast to the standard surge pricing policy, the proposed pricing mechanism results in fewer responsive users to be priced out due to the competition with the unresponsive users despite both pricing mechanisms controlling the size of $q$ very similarly. In all three scenarios in Fig.~\ref{fig:sim_ts}, all three queue sizes stabilize after the burst of price-unresponsive demand. 

The distinction between the use of a vanishing-price versus saturating-price function in our proposed mechanism is that the former trades the rate of convergence to stability for minimizing the number of responsive users being priced out. This is seen by the higher levels of queue size for the price-responsive ($R$) queue and the prolonged presence of the price-unresponsive population ($U$) in the time series of Fig.~\ref{fig:sim_ts}B compared to Fig.~\ref{fig:sim_ts}C.

\begin{figure}[htp]
    \centering
    \includegraphics[width=.6\linewidth]{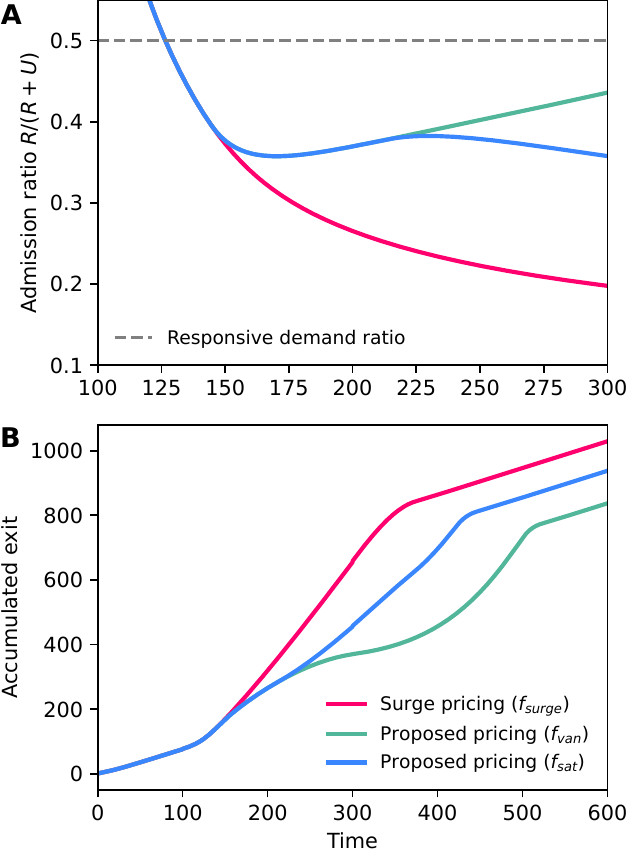}
    \caption{Fairness measurement comparison in simulation scenarios with different pricing mechanisms. The top panel (A) shows Admission ratio of responsive users ($R/(R+U)$) over time, comparing the standard surge pricing (red curve) with the proposed pricing (green and blue curves corresponding to the vanishing-price and saturating-price functions). The dashed line marks the ratio corresponding to admitting the same number of users from each category. The bottom panel (B) compares the three different pricing mechanisms by the total number of price-responsive users that have left the system due to the price, from the start of the simulation to any time $t$.}
    \label{fig:sim_fairness}
\end{figure}

Finally, we consider the effect of the proposed mechanism on fairness. 
Figure~\ref{fig:sim_fairness}A  compares the admission ratio for the price-responsive population for the three pricing functions by depicting the ratio of responsive users in all admitted users to the service queue (from simulations of Fig.~\ref{fig:sim_ts}) over time
using (i) the surge pricing policy (red), 
(ii) the proposed vanishing-price (green) and (iii) the saturating-price (blue) control mechanisms. The increase (resp. decrease) in the admission ratio (of responsive population) using the proposed pricing (resp. surge pricing), corresponds to the return (resp. exit) of the responsive population in Fig.~\ref{fig:sim_ts}B,C (resp. Fig.~\ref{fig:sim_ts}A) at approximately $t=175$ and onward. The difference between surge pricing and the proposed pricing mechanisms---marked by significantly higher responsive admittance under the proposed mechanism---indicates that our design improves fairness when unresponsive demand becomes high enough to undermine the effectiveness of surge pricing in balancing supply and demand. 

In Fig.~\ref{fig:sim_fairness}B, the three pricing functions are compared in terms of the total number of price-responsive users priced out of the system accumulated over time, i.e., $\int_t f(q)R\text{d}t$. The conclusions are similar to that of the previous comparisons, however, \enquote{accumulated exit} provides a more intuitive measure of how the price is affecting the users and governing the competition between the two classes of users. 

\section{Conclusions}
We have considered an alternative to dynamic pricing as a means for controlling access to a shared queue. Our pricing control mechanism is fairer (less socially regressive) when compared to traditional dynamic pricing schemes, while at the same time managing access in a dynamic pricing manner when traffic is homogeneously responsive. Our experimental results illustrate the effectiveness of the proposed schemes. 
Future work will consider generalizations of the pricing schemes considered in the paper.
An intriguing direction building on the present work is the development of schemes to explicitly optimize a pricing function so as to attain desired objectives while provably guaranteeing specification constraints such as preventing given congestion levels or high rates of user dropouts due to pricing. Other research directions to be explored include the development of more refined access and dropout models.

\appendix

\section{Details for the proof of Lemma~\ref{lemma:localstability-2d}}
\label{appendix}

Following suggestions from some readers of earlier versions of this paper, we provide the details of the computations used in the proof of Lemma~\ref{lemma:localstability-2d}.
For the computation of the Jacobian of the vector field
\begin{equation*}
    F(R,q) = \begin{bmatrix}
       F_1(R,q) \\ F_2(R,q) 
    \end{bmatrix}= \begin{bmatrix}
        K_R - f(p(q)) R - \alpha(q) R \\
        \alpha(q)R - \mu(q)
    \end{bmatrix}
\end{equation*}
recall that the Jacobian is defined as
\begin{equation}
\label{eq:appJac}
    DF(R,q) = \begin{bmatrix}
        \frac{\partial}{\partial R} F_1(R,q) & \frac{\partial}{\partial q} F_1(R,q) \\
        \frac{\partial}{\partial R} F_2(R,q) & \frac{\partial}{\partial q} F_2(R,q)
    \end{bmatrix}.
\end{equation}
In the following we use the sum rule of differentiation and the fact that the derivative of a linear function is given by $\frac{\dx}{\dx x} (ax +b) = a$. We also use the notation introduced in Lemma~\ref{lemma:localstability-2d} and write $f'(q) := \frac{\dx}{\dx q} f(q)$ etc.
Using this we get
\begin{align*}
    \frac{\partial F_1}{\partial R} (R,q) &= \frac{\partial}{\partial R} \Big( K_R - f(p(q))R - \alpha(q)R\Big) = - (f(p(q)) + \alpha(q)), \\
    \frac{\partial F_1}{\partial q} (R,q) &= \frac{\partial}{\partial q}\Big( K_R - f(p(q))R - \alpha(q)R\Big) = - ((f\circ p)'(q) + \alpha'(q))R, \\
    \frac{\partial F_2}{\partial R} (R,q) &= \frac{\partial}{\partial R}
    \Big( \alpha(q) R - \mu(q) \Big)= \alpha(q), \\
    \frac{\partial F_2}{\partial R} (R,q) &= \frac{\partial}{\partial q}
    \Big( \alpha(q) R - \mu(q) \Big)= \alpha'(q)R - \mu'(q)
\end{align*}
in all points in which the respective derivatives exist. We have also used that the derivative of the constant $K_R$ is $0$. Arranging these terms in the right order into \eqref{eq:appJac} we obtain \eqref{eq:2dim-linearization}.

By definition, the divergence of the vector field $F$ is the trace of the Jacobian. This results in
\begin{equation*}
    \dive F(R,q) = \frac{\partial F_1}{\partial R} (R,q) + \frac{\partial F_2}{\partial q} (R,q)
\end{equation*}
and using our previous derivations, we obtain in the points of differentiability
\begin{equation}
\label{eq:appdivedef}
  \dive F(R,q) = - (f\circ p)(q) - \alpha(q) + \alpha'(q) R - \mu'(q).  
\end{equation}
Finally, we use the definition of $\mu$ in \eqref{eq:mudef} and the rule for the differentiation of linear functions mentioned earlier to obtain
\begin{equation}
\label{eq:appmuprime}
    \mu'(q) =  \frac{\mu^*}{q_c}, \quad q \in (0,q_c) \quad \mu'(q) = 0 , \quad q > q_c.
\end{equation}
Then \eqref{eq:divergence} follows by inserting \eqref{eq:appmuprime} into \eqref{eq:appdivedef}.

\end{document}